\documentclass[a4paper,11pt]{amsart}

\usepackage{amssymb,enumerate,color}
\usepackage{color,mathrsfs}
\usepackage{latexsym,amsbsy,bm,esint}

\usepackage{placeins}

\input epsf

\renewcommand{\O}{\Omega}
\renewcommand{\a}{\alpha}
\newcommand{\pn}{\phi^\alpha_n}
\newcommand{\pz}{\phi^\alpha_0}
\newcommand{\psn}{\psi^\alpha_n}

\newcommand{\Lak}{L^{\alpha}_n}
\newcommand{\fLn}{\mathfrak{L}^\alpha_n}
\newcommand{\fL}{\mathfrak{L}}
\newcommand{\eln}{\ell^{\alpha}_n}
\newcommand{\el}{\mathscr{L}}

\newcommand {\SC} {{\mathbb C}}

\newcommand {\SL} {{\mathbb L}}

\newcommand {\SR} {{\mathbb R}}
\newcommand {\SRd} {{\SR^d}}
\newcommand {\SRp} {{\SR_+}}

\newcommand {\SZ} {{\mathbb Z}}
\renewcommand {\phi} {{\varphi}}

\newcommand {\al} {{\alpha}}

\newcommand {\e} {{\varepsilon}}
\newcommand {\ga} {{\gamma}}
\newcommand {\Ga} {{\Gamma}}
\newcommand {\la} {{\lambda}}
\newcommand {\lam} {{\lambda}}
\newcommand {\La} {{\Lambda}}

\newcommand{\cL}{\mathcal L}
\newcommand{\cM}{\mathcal M}

\newcommand {\tf} {{\tilde f}}

\newcommand {\bI} {{\bar I}}

\def\supp{\mathop{\rm supp}}

\renewcommand {\l}[1]{\langle{#1}\rangle}
\renewcommand {\la}[1]{\langle{#1}\rangle^{\al+\frac12}}

\newcommand {\lat}[1]{\langle{#1}\rangle^{\al+\frac32}}

\numberwithin{equation}{section}
\newtheorem{theorem}{Theorem}[section]
\newtheorem{lemma}[theorem]{Lemma}

\newtheorem{corollary}[theorem]{Corollary}
\newtheorem{Remark}[theorem]{Remark}
\newtheorem{proposition}[theorem]{Proposition}

\newtheorem{example}[theorem]{Example}

\newtheorem{problem}{Problem}

\newcommand {\Proofof}[1] {\noindent{\bf P{\footnotesize\bf ROOF} of {#1}: } \ }

\newcommand {\Proof} {\noindent{\bf P{\footnotesize\bf ROOF}: } \ }
\newcommand {\ProofEnd} {
             \begin{flushright} \vskip -0.2in $\Box$ \end{flushright}}

\newcommand{\Ba}[1]{\begin{array}{#1}}
\newcommand{\Ea}{\end{array}}
\newcommand{\Be}{\begin{equation}}
\newcommand{\Ee}{\end{equation}}
\newcommand{\Bea}{\begin{eqnarray}}
\newcommand{\Eea}{\end{eqnarray}}
\newcommand{\Beas}{\begin{eqnarray*}}
\newcommand{\Eeas}{\end{eqnarray*}}
\newcommand{\Benu}{\begin{enumerate}}
\newcommand{\Eenu}{\end{enumerate}}
\newcommand{\Bi}{\begin{itemize}}
\newcommand{\Ei}{\end{itemize}}

\newcommand{\BR}{\begin{Remark} \em}
\newcommand{\ER}{\end{Remark}}
\newcommand{\BE}{\begin{example} \em}
\newcommand{\EE}{\end{example}}

\newcommand {\Ds} {\displaystyle}

\newcommand {\T} {\mathtt{T}}

\renewcommand {\ae} {{\quad a.e.\,}}
\newcommand {\sae} {{\,a.e.\,}}
\newcommand {\mand} {{\quad\mbox{and}\quad}}
\renewcommand {\mid} {{\,\,\,\colon\,\,\,}}

\renewcommand{\th}{{\operatorname{th}\,}}

\def\Xint#1{\mathchoice
{\XXint\displaystyle\textstyle{#1}}%
{\XXint\textstyle\scriptstyle{#1}}%
{\XXint\scriptstyle\scriptscriptstyle{#1}}%
{\XXint\scriptscriptstyle\scriptscriptstyle{#1}}%
\!\int}
\def\XXint#1#2#3{{\setbox0=\hbox{$#1{#2#3}{\int}$ }
\vcenter{\hbox{$#2#3$ }}\kern-.6\wd0}}

\def\mint{\Xint-}

\renewcommand{\wp}{w^{-\frac1{p-1}}}
\newcommand{\wpp}{w^{-\frac{p'}{p}}}
\newcommand{\Wp}{W^{-\frac1{p-1}}}
\newcommand{\Wpp}{W^{-\frac{p'}{p}}}

\newcommand{\Ve}{V_\e}

\newcommand{\phiT}{\phi_\T}
\newcommand{\hts}{h^*_{t_0}}

\newcommand{\Ml}{\mathcal{M}^{\rm loc}}

\newcommand{\bline}{{\bigskip

\noindent}}

\newcommand{\sline}{{\smallskip

\noindent}}

\begin{document}

\title[2-weight inequalities for Poisson-Laguerre]{A.e. convergence and 2-weight inequalities for Poisson-Laguerre semigroups}
\author[Garrig\'os, Hartzstein, Signes, Viviani]{G. Garrig\'os, S.Hartzstein, T. Signes, B. Viviani}


\address{Garrig\'os and Signes
\\
Departamento de Matem\'aticas
\\
Universidad de Murcia
\\
30100 Murcia, Spain. \emph{Email}: {\rm \texttt{gustavo.garrigos@um.es, tmsignes@um.es}}} 

\address{Hartzstein and Viviani
\\
 IMAL (UNL-CONICET) y FIQ (Universidad Nacional del Litoral)\\
 CCT CONICET Santa Fe Colectora Ruta Nac. N°168, Paraje El Pozo, 3000 Santa Fe, Argentina.
 \emph{Email}: {\rm \texttt{shartzstein@santafe-conicet.gov.ar}, \texttt{viviani@santafe-conicet.gov.ar}}}

\date{\today}


\begin{abstract}
We find optimal decay estimates for the Poisson kernels
associated with various Laguerre-type operators $L$.  
From these, we solve two problems about the Poisson semigroup $e^{-t\sqrt{L}}$.
First, we find the largest space of initial data $f$ so that $e^{-t\sqrt{L}}f(x)\to f(x)$
at $\sae  x$. Secondly, we characterize the largest class of weights $w$ which admit 2-weight inequalities
of the form $\|\sup_{0<t\leq t_0}|e^{-t\sqrt{L}}f|\,\|_{L^p(v)}\lesssim \|f\|_{L^p(w)}$, for some other weight $v$. 
\end{abstract}

\maketitle

\section{Introduction}\label{intro}

 In this paper we continue the research, started in \cite{HTV,GHSTV}, about
 Poisson integrals associated with certain differential operators $L$, say symmetric and positive in $L^2(\O,\mu)$. Namely,
we are interested in the behavior of \[
u(t,x)=e^{-t\sqrt{L}}f(x) 
\]
as a solution of the elliptic differential equation 
\[ \left\{\Ba{ll}u_{tt}-Lu=0 & \\
 u(0,x)=f(x),& \Ea\right. \mbox{in the half-plane } (0,\infty)\times \O.\]
We shall study two questions, which are closely related among themselves

\sline (i) find the {\bf largest} class of functions $f$ for which $\lim_{t\to0^+}u(t,x)=f(x)$, $\sae x\in \O$;

\sline (ii) establish 2-weight inequalities of the form 
\[
\big\|\sup_{0<t\leq t_0}|u(t,x)|\,\big\|_{L^p(v)}\lesssim \|f\|_{L^p(w)}
\]
for the {\bf largest} class of weights $w$ for which a suitable $v$ with this property exists.

\medskip

In \cite{GHSTV} we considered these questions in full detail when $L$ is the Hermite operator.
In this paper we intend to do the same for the various \emph{Laguerre operators}.
We remark that solving these questions for the {\bf largest} class of weights or initial data is generally not an easy task, requiring {\bf optimal} decay estimates of the Poisson kernels. 
These new estimates have an independent interest and may be useful in other settings; see e.g. recent work by Liu and Sj\"ogren \cite{LiuSjo}. For our purposes, they will provide a.e. convergence for new initial data $f$ compared to Muckenhoupt \cite{muck69} and Stempak \cite{stm94}, and also for larger weight classes compared to those of Nowak in \cite{nowak}.

\

We now state our results, for simplicity in the special case of the 
classical Laguerre operator in $\SRp:=(0,\infty)$ 
\Be
\SL\,=\,-y\,\partial_{yy}\,-\,(\al+1-y)\,\partial_y\, +\,m,  \quad \mbox{ where } \al>-1\mand m\geq 0.
\label{SL}\Ee
We have incorporated a parameter $m\geq0$ which later will allow us to recover other Laguerre-type
operators (after suitable changes of variables).
We shall also consider a slightly more general family of partial differential equations, namely
\Be
\left\{\Ba{ll}u_{tt} \,+\,\frac{1-2\nu}t\,u_t\,=\,\SL u & \\
u(t,0)=f,&\Ea\right. \quad \mbox{ where }\nu>0.
\label{pde}\Ee
These pde's appear in relation with the \emph{fractional operator} $f\mapsto\SL^\nu f$ (see e.g. \cite{StTo}).
  
\medskip

As discussed in \cite{StTo}, a candidate solution to \eqref{pde} is given by the \emph{Poisson-like integral}
\Be
P_tf(x):=\tfrac{t^{2\nu}}{4^\nu\Gamma(\nu)}\,\int_0^\infty
 e^{-\frac{t^2}{4u}}\,\big[e^{-u\SL}f\big](x)\,\frac{du}{u^{1+\nu}},\quad t>0,
\label{poisf}\Ee
which is subordinated to the ``heat'' semigroup $\{e^{-u\SL}\}_{u>0}$. 
Our first goal is to find the most general conditions on a function $f:\SRp\to\SC$ so that 
$P_tf$ is a meaningful solution of \eqref{pde}. These conditions will depend on the
following weight 
\Be
\Phi(y)=\frac {y^\al e^{-y}}{(1+y)^{m}[\ln(e+y)]^{1+\nu}} \quad \mbox{if } m>0,\mand
\Phi(y)=\frac {y^\al e^{-y}}{[\ln(e+y)]^{\nu}}  \mbox{ if } m=0.
\label{Phi}\Ee

\

\begin{theorem}\label{th1SL}
For every $f\in L^1(\Phi)$ the function $u(t,x)=P_tf(x)$ in \eqref{poisf} is defined by an absolutely convergent
integral such that
\Bi\item[(i)] $u(t,x)\in C^\infty((0,\infty)\times\SRp)$ and satisfies the pde \eqref{pde}
\item[(ii)]   $\lim_{t\to0^+}u(t,x)=f(x)$ at a.e. $x\in\SRp$. 
\Ei
Conversely, if a function $f\geq0 $  is such that the integral in \eqref{poisf} is finite for some $(t,x)\in(0,\infty)\times\SRp$, then $f$ must necessarily belong to $L^1(\Phi)$.
\end{theorem}

We recall that Muckenhoupt proved in \cite{muck69} the
pointwise convergence for data $f$ in the smaller space $L^1(y^\al e^{-y}\,dy)$ (in the classical setting, i.e., \eqref{SL} with $m=0$ and \eqref{poisf} with $\nu=1/2$).
Our result, which is sharp on positive $f$, enlarges the space to  \[L^1(y^\al e^{-y}/\sqrt{\log(e+y)}\,dy),\]
and allows to include new initial data such as $f(y)=e^y/[(1+y)^{\al+1}\log(e+y)]$.

\

Our second question concerns more ``quantitative'' 
bounds for the solutions of \eqref{pde}, expressed in terms of the following \emph{local maximal operators}
\Be
\label{max} 
P_{t_0}^*f(x):=\sup_{0<t\leq t_0}\big|P_tf(x)\big|\;, \quad\mbox{with $t_0>0$ fixed.}\Ee  
From our estimates of the Poisson kernels we shall be able to prove that \[
P^*_{t_0}:L^1(\Phi)\to L^s_{\rm loc}\quad \mbox{if $s<1$}, \mand P^*_{t_0}:L^1(\Phi)\cap L^p_{\rm loc}\to L^p_{\rm loc}\quad 
\mbox{if $p>1$.}\]
However, our main interest is to obtain global bounds in $x$, which we shall phrase through the following problem.


\begin{problem} \label{Prob}{\bf A 2-weight problem for the operator $P^*_{t_0}$.}
Given $1<p<\infty$, characterize the class of weights $w(x)>0$ such that
$P^*_{t_0}$ maps $L^p(w)\to L^p(v)$ boundedly, for some other weight $v(x)>0$.
\end{problem}

Our second main result gives a complete answer to Problem \ref{Prob}. For $p\in(1,\infty)$ we define the class of weights 
\Be
D_p(\Phi) \,=\,\Big\{ w(y)>0 \mid 
\big\|w^{-\frac1p}\,\Phi\,\big\|_{L^{p'}(\SRp)}<\infty\Big\}. \label{Dp}
\Ee
Observe that $L^p(w)\subset L^1(\Phi)$ if and only if $w\in D_p(\Phi)$, so in view of Theorem \ref{th1SL}, this is a necessary condition for Problem \ref{Prob}. Our second theorem shows that it is also sufficient.

\begin{theorem}\label{th2SL}
Let $1<p<\infty$ and $t_0>0$ be fixed. Then, for a weight $w(x)>0$ the condition $w\in D_p(\Phi)$ is equivalent to
the existence of some other weight $v(x)>0$  such that \Be
P^*_{t_0}:L^p(w)\to L^p(v) \quad \mbox{boundedly}.\label{bded}\Ee
Moreover, for every $\e>0$, we can choose a weight $v\in D_{p+\e}(\Phi)$ satisfying \eqref{bded}.
\end{theorem}

\

We remark that Problem \ref{Prob} is only a ``one-side'' problem, in contrast with the (more difficult) question of characterizing all pairs of weights $(w,v)$ for which \eqref{bded} holds. One-side problems were considered in the early 80s by Rubio de Francia \cite{RdF1} and Carleson and Jones \cite{CJ} for various classical operators. Here we shall follow the approach by the latter, which has the advantage of giving {\bf explicit} expressions
for the second weight $v(x)$ (see Remark \ref{vPhi} below). This is also a novelty compared to \cite{GHSTV}, where
we used the non-constructive method of Rubio de Francia.

\

We can now briefly describe our approach to the proofs in this paper. Rather than 
working with the operator $\SL$, most of our computations will involve
the ``squared'' Laguerre operator
\Be L=-\partial_{yy} + \Big[y^2 + \frac{\al^2-\tfrac14}{y^2}\Big] \,+\,2\mu.\label{Ls1}\Ee
This has the advantage of resembling the Hermite operator (which is the case $\al=-1/2$), so at some points 
we may use computations from \cite{GHSTV}. There are however various
additional difficulties which are characteristic of the Laguerre setting. The term $1/y^2$
 produces a singularity when $y\to0$ which must be handled separately from the singularity at $y\to\infty$.
This is reflected in the behavior of the Bessel function $I_\al$ which is part of the kernel expression
of $e^{-uL}$. One may also expect additional difficulties when $\al\in(-1,-1/2)$ (cases sometimes avoided in the literature, but that we consider here), related to the fact that such Laguerre functions blow-up when $y\to0$. 

Most of our work will be employed in deriving \emph{precise decay estimates} for the Poisson kernel,
which will lead to the following control of the operator $P^*_{t_0}$\Be
P^*_{t_0}f(x)\,\lesssim\,C(x)\,\Big[\Ml(f\Phi)(x)\,+\,\int_{\SRp}f\Phi\;\Big],
\label{Pcontrol}\Ee
for a reasonably well-behaved $C(x)$ (to be absorbed later as part of the weight $v(x)$).
Here $\Ml$ denotes a \emph{local} Hardy-Littlewood maximal operator in $\SRp$, given by 
\Be
\Ml f(x)\,:=\,\sup_{r>0}\,\frac1{|I(x,r)|}\,\int_{I(x,r)} |f(y)|\,\chi_{\{\frac x2\leq y\leq Mx\}}\,dy
\label{Ma}\Ee
for a suitable $M>1$. We also use the notation 
$I(x,r)=(x-r,x+r)\cap\SRp$.

\medskip

Finally, we remark that the statement of Theorems \ref{th1SL} and \ref{th2SL} remains true when $\SL$ is replaced by any
of the Laguerre-type operators in the table below, provided $\Phi$ in \eqref{Phi} is replaced  by the corresponding function in the table. As in \eqref{Phi}, in the extremal case $\mu=-(\al+1)$, 
the logarithmic term in the denominator of $\Phi$ must be replaced by $[\log(e+y)]^\nu$.


\begin{table}[!ht]
  \centering
  \begin{tabular}{|c|c|c|} \hline
Eigenfunctions & differential operator   & function $\Phi$         \\
\hline
$\Lak $ & $\Ba{l} \\ \left.\Ba{r} \SL=-y\,\partial_{yy}\,-\,(\al+1-y)\,\partial_y\, +\,m\\ m\geq 0\Ea\right\} \\ \\ \Ea$           &
  $\Ds\frac {y^\al e^{-y}}{(1+y)^{m}[\ln(e+y)]^{1+\nu}}$    \\

\hline
 $\pn$     & $\Ba{l} \\ \left.\Ba{r}L=-\partial_{yy} + y^2 + \tfrac1{y^2}\big(\al^2-\tfrac14\big) \,+\,2\mu\\
\mbox{{\small$\mu\geq -(\al+1)$}}\Ea\right\} \\ \\ \Ea$
     &               $\Ds\frac{y^{\alpha +  \frac 12 }\;e^{-y^2 /2}}{(1+y)^{1+\alpha+\mu}\,[\ln(e+y)]^{1+\nu}}$    \\
\hline
$\mathfrak{L}^{\alpha}_n$     & $\Ba{l} \\ \left.\Ba{r}\mathfrak{L}=-y\partial_{yy} -\partial_y+
\tfrac14\Big[y + \tfrac{\al^2}y\Big]\,+\,\tfrac\mu2\\
\mbox{{\small$\mu\geq -(\al+1)$}}\Ea\right\} \\ \\ \Ea$
     &               $\Ds\frac{y^{\frac\alpha2}\;e^{-y/2}}{(1+y)^{\frac{1+\alpha+\mu}2}\,[\ln(e+y)]^{1+\nu}}$    \\   
\hline $\ell^{\alpha}_n$ & $\Ba{l} \\ \left.\Ba{r}\mathscr{L}=-y\partial_{yy} -(\al+1)\partial_y+
\tfrac y4\,+\,\tfrac\mu2\\
\mbox{{\small$\mu\geq -(\al+1)$}}\Ea\right\} \\ \\ \Ea$
     &               $\Ds\frac {y^\al e^{-y/2}}{(1+y)^{\frac{1+\alpha+\mu}2}\,[\ln(e+y)]^{1+\nu}}$    \\
 \hline
$ \psi^{\alpha}_n$     & $\Ba{l} \\ \left.\Ba{r}\La=-\partial_{yy} -\tfrac{2\al+1}y\partial_y\,+ y^2\,+\,2\mu\\
\mbox{{\small$\mu\geq -(\al+1)$}}\Ea\right\} \\ \\ \Ea$
     &               $\Ds\frac{y^{2\alpha + 1}\;e^{-y^2 /2}}{(1+y)^{1+\alpha+\mu}\,[\ln(e+y)]^{1+\nu}}$       \\
\hline

\hline
\end{tabular}

\sline  \caption{Table of $\Phi$-functions for various Laguerre-type operators.}
 \label{table1}
\end{table}


\FloatBarrier

The outline of the paper will be the following. In $\S\ref{sec_heat}$ we consider a version of Theorem \ref{th1SL}
for \emph{heat integrals} $u(t,x)=e^{-tL}f(x)$, which are solutions of the heat equation \[
u_t\,+\,Lu=0\quad \mbox{in }(0,\T)\times\SRp,\quad\mbox{with}\quad u(0,x)=f(x).
\]
Heat integrals are easier to handle, and the explicit expression of the heat kernel, $e^{-tL}(x,y)$,
makes more transparent the behavior we shall later encounter in Poisson kernels. In $\S\ref{secmax}$
we study 2-weight inequalities for the local maximal operator $\Ml$. In $\S\ref{secheatTh2}$ we apply these to prove a version of Theorem \ref{th2SL} for heat integrals. 
In $\S\ref{secPois}$ we take up the study of Poisson integrals, splitting in various subsections the detailed kernel estimates leading to \eqref{Pcontrol}. 
In $\S\ref{proofs}$ we shall give the proof of Theorems \ref{th1SL} and \ref{th2SL} for the operator $L$. Finally, in $\S\ref{transf}$ we show how to transfer the results to the
Laguerre operators in Table \ref{table1}.

\

Throughout the paper $\al\!>\!-1$ is fixed, as are the parameters $\mu,m$ in the differential operators.
The notation $A\lesssim B$ will mean $A\,\leq\, c\,B$, for a constant $c>0$ which may depend on 
$\al,\mu$ and other parameters like $p,M,t_0,\e$, but not on $t,x,y$. If needed, we shall stress the latter dependence by $c(x)$, $c(t,x)$, ... Finally, if $1<p<\infty$ we set $p'=p/(p-1)$.

\section{The simpler model of heat integrals}\label{sec_heat}

In this section $L$ will denote the Laguerre-type operator \Be L=-\partial_{yy} + \Big[y^2 + \frac{\al^2-\tfrac14}{y^2}\Big],\label{Ls2}\Ee
that is, we have set  $\mu=0$ in \eqref{Ls1}\footnote{For heat integrals this implies no loss, since the general case can be factored as $e^{-2\mu t}[e^{-tL}f](x)$.}. The corresponding eigenfunctions $\{\pn\}_{n=0}^\infty$ satisfy \[L\pn = (4n+2\al+2)\pn,\quad n=0,1,2,\ldots\]
and form an orthonormal basis of $L^2(0,\infty)$.
The kernel of the associated heat semigroup $e^{-tL}$, written in terms of the new variable $s=\th t$, has the explicit expression
\Bea e^{-tL}(x,y) & = & \sum_{n=0}^\infty e^{-(4n+2\al+2)t}\pn(x)\pn(y)\nonumber\\
& = & \sqrt{\frac{1-s^2}{2s}}\,\bI_\al\Big(\frac{(1-s^2)xy}{2s}\Big)\,e^{-\frac{(x-y)^2}{4s}}\,e^{-\frac{s(x+y)^2}4}.\label{kerIa}\Eea
Here we have used the convenient notation $\bI_\al(z) = \sqrt{z} e^{-z} I_\al(z)$, so that $\bI_\al(z)  \approx \la{z}$, with  $\l{z}=\min\{z,1\}$.

\subsection{A.e. convergence of heat integrals}

We wish to establish the pointwise convergence of $e^{-tL}f(x)$ with the weakest possible conditions in $f$. For this purpose, the following kernel bound will suffice
\Be
e^{-tL}(x,y)\,\lesssim\,\Big\langle\frac{xy}s\Big\rangle^{\al+\frac12}\,\,\frac{e^{-\frac{(x-y)^2}{4s}}}{\sqrt s}.
\label{hker}\Ee
To produce this bound from \eqref{kerIa} one disregards the last exponential, and uses $1-s^2\leq 1$ when $\al\geq -\frac12$. If $\al\in(-1,-\frac12)$,
note that $\l{\lam z}\geq \lam\l{z}$ for $\lam\leq 1$, so one can leave outside a power $(1-s^2)^{\al+1}\leq 1$.

\begin{theorem}\label{heat1}
Let $\al>-1$ be fixed, and $f$ be such that
\Be
\int_0^\infty |f(y)| e^{-a y^2}\l{y}^{\a+\frac12}\,dy <\infty, \quad \mbox{for some (possibly large) } a>0.
\label{faa}\Ee
Then, 
\[
\lim_{t\to0^+} e^{-tL}f(x) = f(x), \ae x\in\SR_+.\]
\end{theorem}

\Proof  For each fixed $N\geq2$ it suffices to show that 
$\lim_{t\to0^+} e^{-tL}f(x) = f(x)$ for $\sae x\in(1/N,N)$.
We split 
\[
f= f\chi_{\{0<y\leq 2N\}}\,+\,f\chi_{\{y>2N\}}\,=\,f_1\,+f_2.\]
The function $f_1$ has bounded support and belongs to $L^1(y^{\al+\frac12}e^{-\frac{y^2}2}dy)$, so we can apply the results of Muckenhoupt \cite{muck69} (with a suitable change of variables\footnote{See e.g. $\S\S$\ref{subsecfLn} and \ref{subsecLak} below for the explicit change of variables.}, as indicated by Stempak \cite{stm94}) to obtain
\[
\lim_{t\to0^+} e^{-tL}f_1(x) = f_1(x)= f(x), \ae x\in\big[\tfrac1N,N\big].
\]
Next we shall show that, under the hypothesis \eqref{faa},
\[
\lim_{t\to0^+} e^{-tL}f_2(x)=0, \quad\forall x\in\big[\tfrac1N,N\big].
\]
Since $t\to0$, we may assume that $s=\th t \leq s_0$ for some $s_0<\frac1{10}$ (which we shall make precise below). Note that $\tfrac1N\leq x\leq N$ and $y>2N$ imply that 
\[
\Big\langle\frac{xy}{s}\Big\rangle^{\al+\frac12} = 1,\quad \forall s<1.
\]
So, by \eqref{hker}, in this region we have a gaussian bound for the kernel
\[
e^{-tL}(x,y)\, \lesssim\, s^{-\frac12}{e^{-\frac{(x-y)^2}{4s}}}\,\leq \,s^{-\frac12}{e^{-\frac{y^2}{16s}}},\]
using in the last step that $|x-y|\geq y/2$. Choosing $s_0<\frac1{32a}$ (with $a$ as in \eqref{faa}), we see that
for all $y>2N$,
\[
e^{-tL}(x,y)\, \lesssim\, \,s^{-\frac12}{e^{-\frac{y^2}{32s}}}\,e^{-a{y^2}}\,\leq \,s^{-\frac12}{e^{-\frac{N^2}{8s}}}\,e^{-a{y^2}}\]
and therefore
\[
e^{-tL}f_2(x)\, \lesssim\, s^{-\frac12}e^{-\frac{N^2}{8s}}\,\int_{y>2N}|f(y)|\,e^{-a{y^2}}\,dy\longrightarrow 0,\quad\mbox{as } s\to 0^+.\]
\ProofEnd

\subsection{Heat kernel estimates}

The estimates in the proof of Theorem \ref{heat1}, slightly refined in some steps,
lead to the following proposition.

\begin{proposition}\label{heatP1}
Let $\al>-1$. Then,  for every $\ga>1$ there is some $M=M_\ga>1$ such that 
\Be
e^{-tL}(x,y)\, \leq \,
C_{\ga} \left\{\Ba{lll}
\Ds\frac{e^{-\frac{|x-y|^2}{4s}}}{\sqrt s}\;\Big\langle{\frac{xy}s}\Big\rangle^{\al+\frac12}
& {\rm if} &  \frac{x}{2}\leq y\leq Mx\\ & & \\
c(x)\,
\la{y}\,e^{-\frac{y^2}{2\ga\th(2t)}} & {\rm if} &  y<\frac{x}{2} \quad{\rm or}\quad y>Mx\Ea\right.
,\label{etL}\Ee
for all $x,y\in\SRp$ and $s=\th t\in(0,1)$.  Here we can set $c(x)= 1/\l{x}^{\al+\frac32}$. 
\end{proposition}
\Proof
Given $x,y$ and $s$, for simplicity we write $z=\frac{xy}{s}$.
Our estimates below will follow from
\Be
e^{-tL}(x,y)\,\lesssim\,\la{z}\,s^{-\frac12}\,{e^{-\frac{(x-y)^2}{4s}}}\,e^{-\frac{sy^2}4}.
\label{hker2}\Ee
This clearly implies the estimate in the local part $y\in[\frac x2,Mx]$,
so we shall look at the complementary range.
Below we shall ignore the last exponential factor in \eqref{hker2}, and observe that all our estimates
will end up with $e^{-{y^2}/{(4\ga s)}}$. Combining these two one obtains the asserted exponential bound, since\[
e^{-\frac{y^2}{4\ga s}}\,e^{-\frac{sy^2}4}\leq \,e^{-\frac{y^2}{4\ga}(s+\frac1s)}\,=\,e^{-\frac{y^2}{2\ga\th(2t)}},
\] 
as $s+s^{-1}=2/\th(2t)$ when $s=\th t$.

To handle the kernel expression in \eqref{hker2} we need to separate the cases $z\leq 1$ and 
$z\geq 1$. We begin with $z\geq1$. In the region $y>Mx$ we may use $|x-y|\geq (1-\frac1M) y$
to obtain
\[
e^{-tL}(x,y)\, \lesssim\, 
\frac{e^{-(\frac{M-1}M)^2\,\frac{y^2}{4s}}}{\sqrt s}\,\leq \,
c_M\,\frac{e^{-(\frac{M-1}M)^3\,\frac{y^2}{4s}}}{y}\,\leq \,c_M\,{e^{-\frac{y^2}{4s\ga}}}\,
\frac{\la{y}}{\l{x}^{\al+\frac32}},
\]
where in the last step we select $M=M_\ga$ sufficiently large so that $(\frac M{M-1})^3\leq\ga$, and have used  the trivial estimate \Be
\frac1y\,\leq\, \frac{\la{y}}{\l{y}^{\al+\frac32}} \,\leq \, \frac{\la{y}}{\l{x}^{\al+\frac32}},\quad \mbox{if } y\geq x.\label{1y}\Ee
On the other hand, if $y<x/2$ we have $|x-y|\geq x/2$, which leads to
\[
e^{-tL}(x,y)\, \lesssim\, 
\frac{e^{-\frac{(x/2)^2}{4s}}}{\sqrt s}\,\leq \,
c_\ga\,\frac{e^{-\frac{(x/2)^2}{4s\ga}}}{x}\,\leq \,
c_\ga\,\frac{e^{-\frac{y^2}{4\ga s}}}{x}\,.
\]
In the case $\al\in(-1,-\frac12)$ this can be combined with 
\Be
\frac1x\,\leq\, \frac{\la{x}}{\l{x}^{\al+\frac32}} \,\leq \, \frac{\la{y}}{\l{x}^{\al+\frac32}}, \quad\mbox{since }x\geq y.\label{1x}\Ee
If on the contrary $\al\geq -\frac12$, we can insert $1\leq z^{\al+\frac12}$ in the gaussian bound and obtain
\Be
e^{-tL}(x,y)\, \lesssim\, 
\frac{e^{-\frac{(x/2)^2}{4s}}}{\sqrt s}\,\big(\tfrac{xy}{s}\big)^{\al+\frac12}\,= \,
\big(\tfrac{x^2}s\big)^{\al+1}\,{e^{-\frac{(x/2)^2}{4s}}}\,\frac{y^{\al+\frac12}}{x^{\al+\frac32}}\,\leq \,
c_\ga\,{e^{-\frac{y^2}{4s\ga}}}\,\frac{\la{y}}{\l{x}^{\al+\frac32}}\,.\label{etL_aux1}\Ee
This completes the proof of \eqref{etL} when $z\geq1$.

We turn to the case $z\leq 1$, and replace the gaussian bound by\Be
e^{-tL}(x,y)\, \lesssim\, s^{-\frac12}{e^{-\frac{(x-y)^2}{4s}}}\,z^{\al+\frac12}.
\label{etL_aux2}\Ee
This is a better bound when $\al\geq-\frac12$, so some of the previous arguments also lead to \eqref{etL};  namely one can disregard $z$ in the region $y>Mx$, and must keep it when $y<\frac x2$ and argue as in \eqref{etL_aux1}. We are left with the case $\al\in(-1,-\frac12)$, which makes $z^{\al+\frac12}\geq 1$.
In the region $y>Mx$, this can be absorbed by the exponentials since
\Beas 
e^{-tL}(x,y) & \lesssim & \big(\tfrac{xy}s\big)^{\al+\frac12}\,s^{-\frac12}\,
{e^{-(\frac{M-1}M)^2\,\frac{y^2}{4s}}}\,=\, \big(\tfrac{y^2}s\big)^{\al+1}\,
{e^{-(\frac{M-1}M)^2\,\frac{y^2}{4s}}}\,\,\frac{x^{\al+\frac12}}{y^{\al+\frac32}}\\
&\leq &
c_M\,{e^{-(\frac{M-1}M)^3\frac{y^2}{4s}}}\,
\frac{\la{x}\la{y}}{\l{y}^{2\al+2}}\,\leq \, c_M\,{e^{-\frac{y^2}{4s\ga}}}\,
\frac{\la{x}\la{y}}{\l{x}^{2\al+2}}\,,
\Eeas
using in the last step that $y\geq x$. Finally, in the region $y<\frac x2$ the inequalities in \eqref{etL_aux1} remain also valid, so we have completed the proof of Proposition \ref{heatP1}.
\ProofEnd

Proposition \ref{heatP1} can be expressed in terms of the local Hardy-Littlewood maximal function in $\SRp$
\Be
\Ml_M f(x):=\sup_{r>0}\frac1{|I(x,r)|}\int_{I_r(x)} |f(y)|\,\chi_{\{\frac x2<y<Mx\}}\,dy,
\label{MlM}
\Ee
where $I(x,r)=I_r(x)$ denotes the interval $(x-r,x+r)\cap\SRp$.

\begin{corollary}\label{cor3}Let $\al>-1$ and $\ga>1$. Then there is some $M=M_\ga>1$ such that 
\Be
\sup_{0<t\leq t_0} \big|e^{-tL}f(x)\big|\, \lesssim \,
\Ml_M f(x)\,+\,
c(x)\,\int_\SRp|f(y)|\,
\la{y}\,e^{-\frac{y^2}{2\ga\th{(2t_0)}}}\,dy
,\label{Malest}\Ee
for every $x,y,t_0>0$ and  $c(x)= 1/\l{x}^{\al+\frac32}$. 
\end{corollary}
\Proof
We only have to prove the local estimate, and may assume that $\supp f\subset[\frac x2, Mx]$. If $s=\th t\leq x^2$, then
$z=\frac{xy}s\gtrsim 1$ (since $x\approx y$ when $y\in\supp f$), so the gaussian bound of the kernel and
a standard slicing argument easily lead to
\[
 \big|e^{-tL}f(x)\big|\, \lesssim \,\frac{1}{\sqrt s}\,\int_{\SRp}e^{-\frac{(x-y)^2}{4s}}\,|f(y)|\chi_{\{\frac x2\leq y\leq Mx\}}\,dy\,\lesssim\,
\Ml_M f(x)\,.
\]
If $s=\th t\geq x^2$, then using $x\approx y$, 
\[
 \big|e^{-tL}f(x)\big|\, \lesssim \,\frac{x^{2\al+1}}{s^{\al+1}}\,\int_{\frac x2}^{Mx}|f(y)|\,dy\,\lesssim
\, \frac{x^{2\al+1}}{x^{2\al+2}}\,\int_{\frac x2}^{Mx}|f(y)|\,dy\,\lesssim\,\Ml_Mf(x).
\]
\ProofEnd
\BR
An estimate quite similar to \eqref{Malest}, with a slightly worse bound for the exponential inside the integral, was obtained by Chicco-Ruiz and Harboure in \cite[$\S5$]{ChiHar}.
\ER

\section{2-weight inequalities for $\Ml$}\label{secmax}

This section is about the \emph{local maximal operator} in $\SRp$ \[
\Ml f(x):=\sup_{t>0}\,\frac1{|I_t(x)|}\,\int_{I_t(x)} |f(y)|\,\chi_{\{\frac x2<y<Mx\}}\,dy,
\]
where $M>1$ is a fixed parameter. For simplicity, we do not include the subscript $M$ in the 
notation, but the implicit constants appearing below will all depend on $M$.
For the 1-weight theory of this operator we refer to \cite[$\S6$]{NowSte06}.

For each $p\in(1,\infty)$, consider the following family of weights in $\SRp$
\[
D_p^{\rm loc} \, = \,\Big\{W(x)>0\mid \int_J \Wpp\, <\infty,\quad\forall\;J\Subset(0,\infty)\Big\}.  
\]
Associated with $W\in D_p^{\rm loc}$,  we consider a family of weights $\{V_\e\}_{\e>0}$, defined by
\Be
V_\e(x) \,=\,V(x)\,\rho_\e\big[V(x)\big],
\quad\mbox{where }V(x):=\,\Big[\Ml(\Wpp)(x)\Big]^{-\frac{p}{p'}},
\label{Ve}\Ee
and with the notation $\rho_\e(x):=\min\big\{x^\e,x^{-\e}\big\}$.
 Observe that $V_{\e_2} \leq V_{\e_1}\leq V\leq W$ if $\e_1\leq \e_2$. This definition is a slight variant of the one proposed by Carleson and Jones in \cite{CJ}, and leads to the following 2-weight inequalities. 

\begin{theorem}\label{ThVe}
Let $1<p<\infty$ and $W\in D_p^{\rm loc}$. Then for every $\e>0$
\[\Ml:L^p(W)\to L^p(\Ve)\quad \mbox{boundedly,}\]
where $\Ve$ is defined as in \eqref{Ve}.
\end{theorem}
\Proof
The argument of the proof is due to Carleson and Jones \cite{CJ} (see also a recent application in \cite[Prop. 4.2]{gar}). For completeness, we sketch the modifications required for the local operator $\Ml$.
Call $E_n=\{x\in\SRp\mid \Ml(\Wp)(x)<2^n\}$, $n\in\SZ$, and define the operators
\Be
T_ng(x):=\chi_{E_n}\,\Ml(\Wp g)(x).\Ee
Note that $T_n:L^1(\Wp)\to L^{1,\infty}$, with a uniform bound  in $n$, since
\Be
\Big|\Big\{T_ng(x)>R\Big\}\Big|\leq \Big|\Big\{\mathcal{M}(\Wp g)(x)>R\Big\}\Big|\leq \frac{c_0}R\int_{\SRp}\Wp|g|,
\label{Tn0}\Ee
using in the last step the weak-1  boundedness of the Hardy-Littlewood maximal operator.
Similarly, $T_n:L^\infty(\Wp)\to L^\infty$ with $\|T_n\|\leq 2^n$, since
\Be
\big\|T_ng\big\|_\infty=\sup_{x\in E_n}\big|\Ml(\Wp g)(x)\big|\,\leq 2^n\,\|g\|_\infty.\label{Tninf}\Ee
Thus, by the Marcinkiewicz interpolation theorem we obtain 
\Be
\int_{E_n}|T_n(g)|^p\leq\, c_0\,2^{\frac{np}{p'}}\,\int_{\SRd}|g|^p\,\Wp,\quad n\in\SZ.\Ee
Setting $g=fW^{\frac1{p-1}}$ in the above inequality, this is the same as
\Be
\int_{E_n}|\Ml(f)|^p\leq\, c_0\,2^{\frac{np}{p'}}\,\int_{\SRp}|f|^p\,W,\quad n\in\SZ.\label{Ml_ipol}\Ee
Now, modulo null sets $\SRp=\cup_{n\in\SZ}\big[E_n\setminus E_{n-1}\big]$ (since $0<\Ml(\Wp)(x)<\infty$ at $\sae x$), and we have\[
\Ve(x)\approx 2^{-\frac{np}{p'}} 2^{-\frac{\e|n|p}{p'}},\mbox{ if }x\in E_n\setminus E_{n-1}.
\]
Therefore, we obtain
\Beas \int_{\SRp}|\Ml f|^p\,\Ve & \lesssim &\,
\sum_{n\in\SZ}2^{-\frac{np}{p'}}2^{-\frac{\e|n|p}{p'}}\int_{E_n}|\Ml f|^p\\  
\mbox{{\tiny (by \eqref{Ml_ipol})}}& \lesssim & \Big(
\sum_{n\in\SZ}2^{-\frac{\e|n|p}{p'}}\Big)\,\int_{\SRd}|f|^pW,\Eeas  
as we wished to show.
\ProofEnd

The weight $\Ve$ inherits some of the integrability behavior of $W$ if $\e$ is sufficiently small.
 To state this we first define the subclasses
\Beas
D_p^{0}(\beta) & = & \Big\{W\in D_p^{\rm loc}\mid \int_0^1 \Wpp(y)\,\l{y}^{\beta p'}\,dy <\infty\Big\},\quad \mbox{for }\beta>-1, \\
D_p^{\rm exp}(a) & = & \Big\{W\in D_p^{\rm loc}\mid \int_1^\infty \Wpp(y)\,
e^{-ay^2p'}\,dy <\infty\Big\},\quad \mbox{for }a>0. \Eeas

\begin{proposition}\label{P5.2}
Let $1<p<\infty$ and $W\in D_p^{\rm loc}$. Then, for each $\e>0$, the weight defined in \eqref{Ve} satisfies
 $\Ve\in D^{\rm loc}_q$, for all $\;q>p+{\e p}/{p'}$. Moreover, we additionally have

\medskip
\Benu
\item[(i)]  $W\in D^0_p(\beta)$ implies $\Ve\in D^0_q(\beta)$, 
provided $q>p+\e\,\frac{p}{p'}\,\frac{|1+\beta p'|}{1+\beta}$\\
\item[(ii)] $W\in D^{\rm exp}_p(a)$ implies $\Ve\in D^{\rm exp}_q(b)$, provided
 $q>p(1+\e) M^2 a/b$.\Eenu
\end{proposition}
\Proof Observe that \Be V_\e(x)^{-\frac{q'}q}=\max_{\pm}\big[\Ml(\Wp)(x)\big]^{\frac{p-1}{q-1}(1\pm\e)}.\label{Veq}\Ee
The assumption $q>p+\frac{\e p}{p'}$ implies that $s=\frac{(p-1)(1+\e)}{q-1}<1$. Then, given $J=[a,b]\Subset\SRp$,
\Beas\int_J\big[\Ml(\Wp)\big]^s & \lesssim & |J|^{1-s}\,\big\|\mathcal{M}(\Wp\chi_{J^*})\big\|_{L^{1,\infty}}^s\\ & \lesssim &   c_J\,\Big(\int_{J^*}\Wp\Big)^s<\infty, 
\Eeas
with $J^*=[a/2,Mb]\Subset(0,\infty)$. The same applies if we set $s=\frac{(p-1)(1-\e)}{q-1}$ (which is also $<1$),
so we deduce from \eqref{Veq} that $\int_J\Ve^{-\frac1{q-1}}<\infty$.

We next prove (i), and as before set $s=\frac{(p-1)(1+\e)}{q-1}<1$.
Then, denoting $I_j=[2^{-j-1},2^{-j}]$, we have
\Beas
\int_0^1\big[\Ml(\Wp)\big]^s\,\l{y}^{\beta q'}\,dy
& \lesssim &  \sum_{j=0}^\infty2^{-j\beta q'}\,|I_j|^{1-s}\,
\big\|\mathcal{M}(\Wp\chi_{I_j^*})\big\|_{L^{1,\infty}}^s\\
& \lesssim &   \sum_{j=0}^\infty2^{-j\beta q'}2^{-j(1-s)}
\,\Big(\int_{2^{-j-2}}^{M2^{-j}}\Wp\Big)^s\\
& \hskip -4cm\lesssim &  \hskip -2cm \sum_{j=0}^\infty2^{-j[\beta (q'-p's)+1-s]}
\,\Big(\int_{0}^{M}\Wp(y)\,\l{y}^{\beta p'}\,dy\Big)^s.
\Eeas
This is a finite expression provided\[
\beta (q'-p's)+1-s>0\;,
\]
and using the value of $s=\frac{(p-1)(1+\e)}{q-1}$ and solving for $q$ this is equivalent to
\[
q>p+\frac{\e p(1+\beta p')}{p'(1+\beta)}.\]
In order to have $\int_0^1\Ve^{-\frac1{q-1}}\,\l{y}^{\beta q'}\,dy<\infty$ the previous relation must also hold with $\e$ replaced by $-\e$, so a sufficient condition is 
\[q>p+\frac{\e p\,|1+\beta p'|}{p'(1+\beta)},\] as we wished to show.

We now prove (ii). Let $\ga>1$ (to be precised later), and as before set 
$I_j=[\ga^j,\ga^{j+1}]$ and $s=\frac{(p-1)(1+\e)}{q-1}<1$. Then
\Beas
\int_1^\infty\big[\Ml(\Wp)\big]^s\,e^{-by^2q'}\,dy  
& \lesssim &  \sum_{j=0}^\infty\,e^{-b\ga^{2j}q'}\,\ga^{(1-s)j}\,
\big\|\mathcal{M}(\Wp\chi_{I_j^*})\big\|_{L^{1,\infty}}^s\\
& \hskip -4cm\lesssim &  \hskip -2cm \sum_{j=0}^\infty e^{-b\ga^{2j}q'}\,\ga^{(1-s)j}\,
\,\Big(\int_{\ga^{j}/2}^{M\ga^{j+1}}\Wp\,dy\Big)^s\,\\
& \hskip -4cm\leq &  \hskip -2cm \sum_{j=0}^\infty \,\ga^{(1-s)j}\,e^{-\ga^{2j}[bq'-p'aM^2\ga^{2}s]}
\,\Big(\int_{1/2}^{\infty}\Wp\,e^{-p'ay^2}\,dy\Big)^s.
\Eeas
This is now a finite expression provided
\[
bq'> p'aM^2\ga^{2}s,
\]
which using the value of $s$ and solving for $q$ gives
\[
q>p(1+\e)M^2\ga^{2}a/b.\] Clearly, we can choose a $\ga>1$ with this property under the assumption
\[
q>p(1+\e)M^2a/b.\]
Since this also implies the validity of the estimates with  $\e$ replaced by $-\e$,
we may conclude that $\Ve\in D^{\rm exp}_q(b)$, as desired.
\ProofEnd

\section{2-weight inequalities for local maximal heat operators}\label{secheatTh2}

Let $L$ be as in \eqref{Ls2}, and for each $t_0>0$, consider \[
h^*_{t_0}f(x):=\sup_{0<t\leq t_0} \big|e^{-tL}f(x)\big|.\]
Given any $\T>t_0$, this operator is well defined over functions $f\in L^1(\phi_\T)$, where
\[
\phi_\T(y)=\la{y}\,e^{-\frac{y^2}{2\th(2\T)}}.
\]
We wish to study 2-weight inequalities for $h^*_{t_0}$ over subspaces $L^p(w)\subset L^1(\phi_\T)$.
By duality, the class of weights for which such inclusion holds is given by
\[
D_p(\phi_\T):=\Big\{w>0\mid \big\|w^{-\frac1p}\phi_\T\big\|_{p'}<\infty\Big\}.\]
Here we show that for all such weights the operator $h^*_{t_0}$ satisfies a 2-weight inequality.

\begin{theorem}\label{th2heat}
Let $\T>t_0>0$ and $1<p<\infty$. Then, for every $w\in D_p(\phi_\T)$ there exists another weight $v(x)>0$
such that
\[
h^*_{t_0}:L^p(w)\to L^p(v),\quad \mbox{boundedly}.
\]
Moreover, if $q>p$ and $t_0$ is sufficiently small, then we can select $v\in D_q(\phiT)$.\end{theorem}
\BR The second weight $v(x)$ will be constructed explicitly; see \eqref{v2}, \eqref{v1} and \eqref{Vexplicit} below. Observe that $v$ depends on $\al, p, t_0, \T$ and of course $w$.\ER

\Proofof{Theorem \ref{th2heat}}
The crucial estimate was already given in Corollary \ref{cor3}. We shall use it with the parameter $\ga=\th(2\T)/\th(2t_0)>1$,
which produces a suitable $M=M_\ga>1$ such that
\Be
\hts f(x)\, \lesssim \,
\Ml_M f(x)\,+\,
c(x)\,\int_0^\infty|f(y)|\,\phiT(y)\,dy.
\label{max1}\Ee
The last integral is bounded by $\|f\|_{L^p(w)}\|w^{-\frac1p}\phiT\|_{p'}$, so the second term
will be fine for any weight $v(x)$ such that $c(x)=1/\lat{x}\in L^p(v)$.
For instance we may take
\Be
v_2(x)=\frac{\l{x}^{(\al+\frac32)p-1}}{[\log(e/\l{x})]^2\,(1+x)}
\label{v2}\Ee
which clearly satisfies \[
\int_0^\infty |c(x)|^p\,v_2(x)\,dx\,=\,\int_0^\infty\frac{dx}{\l{x}[\log(e/\l{x})]^2\,(1+x)^p}\,<\,\infty.\]
Further, we claim that $v_2\in D_q(\phiT)$ iff $q>p$. Indeed, since $\phiT(x)$ decays exponentially, it suffices to 
check the integrability for $x$ near 0. Writing $\beta=\al+\frac12$ so that\[
\phiT(x)\approx\l{x}^\beta\mand v_2(x)\approx \frac{\l{x}^{(\beta+1)p-1}}{[\log(e/\l{x})]^2}
\] we easily see that \Be
\int_0^1v_2(x)^{-\frac {q'}q}\,\l{x}^{\beta q'}\,dx\,\approx\,\int_0^1\frac{[\log(e/x)]^{2q'/q}}
{x^{1-q'(\beta+1)(1-\frac pq)}}\,dx\,<\,\infty.\label{Dqv1}\Ee

For the first term in \eqref{max1} we shall use the results in $\S\ref{secmax}$. We first note that
\Be
w\in D_p(\phiT) \quad \Longleftrightarrow\quad w\in D^0_p(\al+\tfrac12)\cap D^{\rm exp}_p(a),
\quad\mbox{with }a=1/(2\,\th 2\T),\label{wW}\Ee 
where the weight classes $D^0_p(\beta)$ and $D^{\rm exp}_p(a)$ were defined just before Proposition \ref{P5.2}.
Then, for every $\e>0$ Theorem \ref{ThVe} gives
\[
\big\|\Ml_M f\big\|_{L^p(v_{1,\e})}\,\lesssim\,\big\|f\big\|_{L^p(w)},
\] provided 
\Be v_{1,\e}(x)=\mathscr{V}(x)\rho_{\e}\Big({\mathscr{V}(x)}\Big),\quad\mbox{where }\mathscr{V}(x)=\Big[\Ml_M\big(w^{\frac1{1-p}}\big)(x)\Big]^{1-p}
\label{v1}\Ee
(or $v_{1,\e}=V_{\e}$ in the notation of \eqref{Ve}). Hence setting \Be v(x)=\min\{v_{1,\e}(x),v_2(x)\}\label{Vexplicit}\Ee
with $v_{1,\e}$ and $v_2$ defined as in \eqref{v1} and \eqref{v2}, we 
have proved that $\hts:L^p(w)\to L^p(v)$.

\medskip

 It remains to verify the last statement in Theorem \ref{th2heat}.
We already know that, for every $q>p$, we have  
$v_2\in D_q(\phiT)$. Concerning $v_{1,\e}$, from the equivalence in \eqref{wW}
it suffices to prove that $\Ve\in D^0_{q}(\al+\frac12)\cap D^{\rm exp}_{q}(a)$ for a sufficiently small $\e$.
The first assertion is immediate from (i) in Proposition \ref{P5.2}. However, (ii) in the same proposition only 
gives $\Ve\in D^{\rm exp}_{\rho}(a)$ if $\rho>p(1+\e)M^2$, where $M=M_\ga$ is the parameter obtained in Proposition \ref{heatP1} by the rule $(\frac M{M-1})^3=\ga=\th(2\T)/\th(2t_0)$. If we allow both $\e$ and $t_0$ be sufficiently small (so that $M$ becomes close enough to 1), then we can set $\rho=q$, and hence conclude that $v_{1,\e}\in D_q(\phiT)$ as desired.
\ProofEnd

\section{Poisson kernel estimates}\label{secPois}

In this section we fix $\al>-1$ and $\mu\geq-(\al+1)$, and
consider the Laguerre-type operator \Be L=-\partial_{yy} + \Big[y^2 + \frac{\al^2-\tfrac14}{y^2}\Big] \,+\,2\mu,\label{Lphi}\Ee
whose eigenfunctions $\{\pn\}_{n=0}^\infty$ form an orthonormal basis of $L^2(0,\infty)$, 
and satisfy \[ L\pn = (4n+2(\al+1+\mu))\pn,\quad n=0,1,2,\ldots\]
They can be expressed in terms of the (normalized) Laguerre polynomials $\Lak$ by\Be
\pn(y)\,=\,\sqrt{2}\,y^{\al+\frac12}\,e^{-\frac{y^2}2}\,\Lak(y^2),
\label{pn}\Ee
although we shall not use this formula here. 
The kernel of the associated heat semigroup, $e^{-tL}$, can be written explicitly in various forms 
\Bea e^{-tL}(x,y) & = & \sum_{n=0}^\infty e^{-[4n+2(\al+1+\mu)]t}\pn(x)\pn(y)\nonumber\\
\mbox{{\footnotesize ($r=e^{-2t}$)}}& = & r^\mu\,\sqrt{\tfrac{2r}{1-r^2}}\,\bI_\al\Big(\frac{2rxy}{1-r^2}\Big)\,
e^{-\frac{(x-ry)^2}{1-r^2}}\,e^{\frac{x^2-y^2}2}\label{hr}\\
\mbox{{\footnotesize ($s=\th t$)}}& = & \Big(\tfrac{1-s}{1+s}\Big)^\mu\,\sqrt{\tfrac{1-s^2}{2s}}\,\bI_\al\Big(\frac{(1-s^2)xy}{2s}\Big)\,e^{-\frac{(x-y)^2}{4s}}\,e^{-\frac{s(x+y)^2}4}\label{hs}\Eea
where as before we have set $\bI_\al(z) = \sqrt{z} e^{-z} I_\al(z)$.
Thus, using the notation $\l{z}=\min\{z,1\}$, we shall have
$\bI_\al(z)  \approx \la{z}$.
Both expressions of the heat kernel will be useful in our later computations. For instance, \eqref{hr} is good when $r\approx 0$,
as it isolates correctly the decaying factor $\la{y}e^{-y^2/2}$. On the other hand, \eqref{hs} will be useful when $s\approx0$ (hence $r\approx 1$), since it makes transparent the gaussian behavior of the singularity $s^{-\frac12}e^{-\frac{(x-y)^2}{4s}}$.

\

Using the subordination formula in \eqref{poisf}, the Poisson kernel associated with $L$ becomes
\Be
P_t(x,y):=\tfrac{t^{2\nu}}{4^\nu\Gamma(\nu)}\,\int_0^\infty
 e^{-\frac{t^2}{4u}}\,\big[e^{-uL}(x,y)\big]\,\frac{du}{u^{1+\nu}},\quad t>0.
\label{Ptxy}\Ee
Changing variables $r=e^{-2u}$ (i.e., $u=\frac12\ln\frac1r$) one sees that
\[
P_t(x,y)\,\approx\,t^{2\nu}\,e^{\frac{x^2-y^2}2}\;\int_0^1
 e^{-\frac{t^2}{2\ln\frac1r}}\, r^{\mu+\frac12}\,
\frac{e^{-\frac{(x-ry)^2}{1-r^2}}}{\sqrt{1-r}}\,\frac{\la{\frac{rxy}{1-r^2}}}{(\ln\frac1r)^{1+\nu}}\,
\frac{dr}{r}\,.
\]
We shall consider two regions of integration according to the behavior of $z:=\Ds\frac{rxy}{1-r^2}$.
The regions will be separated by the number\[
r_0(xy)\,=\,\left\{\Ba{ll}\frac1{2xy} & ,\mbox{ if }xy\geq1\\1-\frac{xy}2 &, \mbox{ if }xy\leq1.\Ea\right.\]
Indeed, it is elementary to check that
\Benu
\item If $\;0<r\leq r_0(xy)\;$ then $\;z\leq1$.
\item If $\;r_0(xy)\leq r<1\;$ then $\;z\geq1/2$.
\Eenu
Thus we can write 
\Beas
P_t(x,y) & \approx & t^{2\nu}\,e^{\frac{x^2-y^2}2}\;\Big[\int_0^{r_0(xy)}
 \,\cdots\,\Big(\frac{rxy}{1-r^2}\Big)^{\al+\frac12}\,
\frac{dr}{r}\,+\,\int_{r_0(xy)}^1
 \,\cdots\,\frac{dr}{r}\;\Big]\\
& = & B_t(x,y)\,+\,A_t(x,y).\Eeas
The next two propositions summarize the estimates we shall need to handle these kernels.
We shall make extensive use of the function
\Be
\Phi(y)\,:=\,\frac{\la{y}\,e^{-y^2/2}}{(1+y)^{\mu+\frac12}\,[\log(y+e)]^{1+\nu}},\label{phiL}\Ee
with the agreement that in the extreme case $\mu=-(\al+1)$ the log in the denominator is just $[\log(y+e)]^{\nu}$.
The first result gives, for fixed $t$ and $x$, the \emph{optimal decay} of $y\mapsto P_t(x,y)$ in terms of the function $\Phi(y)$.

\begin{proposition}\label{P3.1}
Given $t,x>0$, there exist $c_1(t,x)>0$ and $c_2(t,x)>0$ such
that \Be c_1(t,x)\, \Phi(y) \,\leq\,P_t(x,y)\,\leq\,
c_2(t,x)\, \Phi(y) \;,\quad \forall\;y\in\SRp. \label{pk1}\Ee
\end{proposition}

The second result is a refinement of the upper bound in \eqref{pk1} with a few advantages: it is uniform in the variable $t$,  it isolates in the ``local part'' the singularities of the kernel $P_t(x,y)$, and finally provides ``reasonable'' bounds for the constant's dependence on $x$.

\begin{proposition}
\label{P3.2} There exists $M>1$ such that the following holds for all $t,x,y>0$

\Be
P_t(x,y) \, \lesssim \;
C_1(x)\,\frac{t^{2\nu}\,e^{-{y^2}/2}}{\big(t+|x-y|\big)^{1+2\nu}}
\,\chi_{\big\{\frac x2<y<
Mx\big\}}\;+\; C_2(x)\,(t\vee1)^{2\nu}\,\Phi(y),\label{pk2}
\Ee
where $C_1(x) \,=\,(1+x)^{2\nu}e^{\frac{x^2}2}$
and $C_2(x)= [\log(e+x)]^{1+\nu}\,{(1+x)^{|\mu+\frac12|}\,e^{\frac{x^2}2}}/{\lat{x}}$.
\end{proposition}

If we consider, for fixed $M>1$, the \emph{local maximal function} in $\SRp$
\Be
\Ml_M f(x):=\sup_{t>0}\,\frac1{|I_t(x)|}\,\int_{I_t(x)} |f(y)|\,\chi_{\{\frac x2<y<Mx\}}\,dy,
\label{Ml}
\Ee
then we may express \eqref{pk2} as follows.

\begin{corollary}\label{C3.3}
Let $t_0>0$ be fixed. Then there is some $M>1$ such that 
\Be
 P^*_{t_0}f(x) \,\lesssim \,
C_1(x)\,\Ml_M \big(fe^{-\frac{y^2}2}\big)(x)\,+\,C_2(x)\,\big\|f\big\|_{L^1(\Phi)}
,\quad x\in\SRp,\label{P**}
\Ee
with $C_1(x)$ and $C_2(x)$ as in Proposition \ref{P3.2}.
 \end{corollary}

This is the key estimate from which we shall deduce the
theorems claimed in $\S1$. The interested reader may wish to skip 
the technical proofs of the propositions in the next subsections and pass
directly to $\S\ref{proofs}$ for the proof of the theorems.

\subsection{Estimates from below for $B_t(x,y)$}
Recall that
\Be
B_t(x,y)\,\approx\,t^{2\nu}\,(xy)^{\al+\frac12}\,e^{\frac{x^2-y^2}2}\;\int_0^{r_0(xy)}
\frac{r^{\al+\mu+1}}{(1-r)^{\al+1}\,(\ln\frac1r)^{1+\nu}}\,
e^{-\frac{t^2}{2\ln\frac1r}}\,e^{-\frac{(x-ry)^2}{1-r^2}}\,
\frac{dr}{r}\,.
\label{Btxy}
\Ee
The lower bound in Proposition \ref{P3.1} will be obtained by just looking at this integral.

\begin{lemma}
\label{L1}
For fixed $t,x>0$ it holds
\[
B_t(x,y)\,\geq\,c_1(t,x)\,\Phi(y),\quad\forall\;y\in\SRp, 
\]
for a suitable function $c_1(t,x)>0$.
\end{lemma}
\Proof
We first look at $y<x\wedge \frac1x$. Then $xy<1$, and hence $r_0(xy)>1/2$. So, we can estimate
$B_t(x,y)$ by an integral over $0<r<\frac12$, which disregarding irrelevant terms becomes
\[
B_t(x,y)\gtrsim\,t^{2\nu}\,(xy)^{\al+\frac12}\,e^{\frac{x^2-y^2}2}\,\int_0^{\frac12}\frac{r^{\al+\mu+1}}{(\ln\frac1r)^{1+\nu}}\,
e^{-\frac{t^2}{2\ln\frac1r}}\, e^{-\frac{(x-ry)^2}{1-r^2}}\,\frac{dr}r.\]
We can get rid of the first two exponentials using\[
e^{\frac{x^2-y^2}2}\geq 1\quad\mbox{(since $y\leq x$)}\mand 
e^{-\frac{t^2}{2\ln\frac1r}}\geq e^{-\frac{t^2}{2\ln2}} \quad\mbox{(since $r\leq \frac12$)}.
\]
For the last exponential notice that $0<x-ry\leq x$, and hence 
$e^{-\frac{(x-ry)^2}{1-r^2}}\geq e^{-\frac43x^2}$. This leaves a convergent integral in $r$, so we conclude that
\[
B_t(x,y)\gtrsim\,c_1(t,x)\,\la{y},\]
with $c_1(t,x)=t^{2\nu}x^{\al+\frac12}e^{-\frac{t^2}{2\ln2}}e^{-\frac43x^2}$.
Notice that $y\leq 1$ in this range, so we find the required expression for $\Phi(y)$.
\

Suppose now that $y\geq x\vee \frac1x$. Then $xy\geq1$, and hence $r_0(xy)=\frac1{2xy}\leq 1/2$. Arguing
as before we can estimate
$B_t(x,y)$ by
\[
B_t(x,y)\gtrsim\,t^{2\nu}\,(xy)^{\al+\frac12}\,e^{-\frac{y^2}2}\,e^{-\frac{t^2}{2\ln2}}\,
\int_0^{\frac1{2xy}}\frac{r^{\al+\mu+1}}{(\ln\frac1r)^{1+\nu}}\,\, e^{-\frac{(x-ry)^2}{1-r^2}}\,\frac{dr}r.\]
This time we get rid of the exponential inside the integral using\[
|x-ry|\leq x+ry\leq  x+\frac1{2x} \quad\mbox{(since $r\leq \frac1{2xy}$)},
\]
which implies $e^{-\frac{(x-ry)^2}{1-r^2}}\geq e^{-\frac43(x+\frac1x)^2}$. We can easily compute the integral
\[
\int_0^{\frac1{2xy}}\frac{r^{\al+\mu+1}}{(\ln\frac1r)^{1+\nu}}\,\frac{dr}r\;\approx\;\frac1{(xy)^{\al+\mu+1}\,[\log 2xy]^{1+\nu}},\quad \mbox{if  {\small $\al+\mu+1>0$},}
\]
with the right hand side becoming $1/[\log 2xy]^{\nu}$ in the extreme case {\small $\al+1+\mu=0$}.
Since $y>\max\{x,1\}$, note that 
\[
\log(2xy)\,\leq \,\log(2y^2)\,\approx\log(y+e).\]
Thus, combining all the previous estimates we conclude that
\[
B_t(x,y)\gtrsim\,c_1(t,x)\,\frac{e^{-\frac{y^2}2}}{y^{\mu+\frac12}\,[\log(y+e)]^{1+\nu}},\]
which, since $y\geq 1$, is the required expression for $\Phi(y)$ (with the usual agreement when $\mu+\al+1=0$).
In this part we have set  $c_1(t,x)=t^{2\nu}e^{-\frac{t^2}{2\ln2}}x^{-(\mu+\frac12)}\,e^{-\frac43(x+\frac1x)^2}$.

Finally, since the function $y\mapsto P_t(x,y)/\Phi(y)$
is continuous and positive, it is also bounded from below by some $c_1(t,x)$ when
$y$ belongs to the compact set $[x\wedge\frac1x,x\vee\frac1x]$.
\ProofEnd

\subsection{Estimates from above for $B_t(x,y)$}

The next lemma, combined with the previous one, shows that for fixed $t$ and $x$, the function $B_t(x,y)$
essentially behaves like $\Phi(y)$.

\begin{lemma}\label{L2}
For fixed $t,x>0$ it holds
\Be
B_t(x,y)\,\lesssim\,c(x)\,\max\{t^{2\nu},1\}\,\Phi(y),\quad\forall\;y\in\SRp, 
\label{B_L2}\Ee
with $c(x)=1/\lat{x}$.
\end{lemma}
\Proof
We first notice that the two exponential terms in \eqref{hr} can be written as
\Be
e^{-\frac{(x-ry)^2}{1-r^2}}\,e^{\frac{x^2-y^2}2}\,=\,e^{-\frac{1+r^2}{1-r^2}\frac{x^2+y^2}2}\,e^{\frac{2rxy}{1-r^2}}
\,\lesssim\,e^{-\frac{x^2+y^2}2},
\label{aux1_L2}\Ee
since $\frac{1+r^2}{1-r^2}\geq 1$ and in the region of integration of $B_t(x,y)$ the exponent $z=\frac{2rxy}{1-r^2}\lesssim1$. We now separate cases.

\sline (i) \emph{Case $xy\geq1$:} then $r_0(xy)=\frac1{2xy}\leq \frac12$ and
\Be
B_t(x,y)\,\lesssim\,t^{2\nu}\,(xy)^{\al+\frac12}\, e^{-\frac{x^2+y^2}2}\,\int_0^{\frac1{2xy}}
\frac{r^{\al+\mu+1}}{(\ln\frac1r)^{\nu+1}}\,\frac{dr}r.\label{aux4_L2}\Ee
The last integral is approximately given by
\[
\int_0^{\frac1{2xy}}
\frac{r^{\al+\mu+1}}{(\ln\frac1r)^{\nu+1}}\,\frac{dr}r\,\approx \,\Big(\frac1{xy}\Big)^{\al+\mu+1}\,\frac1{[\log(2xy)]^{1+\nu}}
\]
(with the usual convention when $\al+\mu+1=0$ of reducing the log by one power). 
This is a good estimate if we assume that $x\geq 1/2$, since we may use 
\[
\log(2xy)\,\gtrsim \, \log(y\vee 2)\,\approx\,\log(y+e),
\]
and overall obtain
\[
B_t(x,y)\,\lesssim\,t^{2\nu}\,x^{-(\mu+\frac12)}\,e^{-\frac{x^2}2}\,
\frac{y^{\al+\frac12}e^{-y^2/2}}{y^{\al+\mu+1}[\log(y+e)]^{1+\nu}}\,\lesssim\,t^{2\nu}\,\Phi(y).
\]
When $x\leq1/2$, we need a refinement to obtain the $c(x)$ in the statement of the lemma.
We split the integral defining $B_t(x,y)$ as\Be
B_t(x,y)\,= \,\int_0^{\frac{2x}y}\cdots \;+\;\int_{\frac{2x}y}^{r_0(xy)}\cdots\;=I\,+\,II,
\label{aux3_L2}\Ee
noticing that the partition point $\frac{2x}y\leq r_0(xy)=\frac1{2xy}$.
Since $x\leq \frac12$ and $xy\geq 1$ we also have $y\geq 2$. Now, the first integral can be bound as above by 
\Beas
I & \lesssim  & t^{2\nu}\,(xy)^{\al+\frac12}\, e^{-\frac{x^2+y^2}2}\,\Big(\frac xy\Big)^{\al+\mu+1}\,\frac1{[\log(\frac y{2x})]^{1+\nu}}\\
 & \lesssim & 
t^{2\nu}\,\la{x}\,x^{\al+\mu+1}\,e^{-\frac{x^2}2}\,
\frac{e^{-y^2/2}}{y^{\mu+\frac12}[\log y]^{1+\nu}}\,\lesssim\,t^{2\nu}\,\la{x}\,\Phi(y).
\Eeas
since in this range $y\geq 2$. This implies the stated estimate because $\la{x}\leq 1/\lat{x}=c(x)$. 
To handle $II$ we need a different bound for the exponentials in \eqref{aux1_L2}, noticing that
\Be
r>\tfrac{2x}y\quad\Rightarrow\quad|x-ry|= ry-x\,\geq\,\tfrac{ry}2\quad\Rightarrow\quad
e^{-\frac{(x-ry)^2}{1-r^2}}\,\leq\,e^{-\frac{r^2y^2}4}.
\label{aux6_L2}\Ee
Thus
\Be
II \; \lesssim \;
t^{2\nu}\,(xy)^{\al+\frac12}\, e^{\frac{x^2-y^2}2}\,
\int_{\frac{2x}y}^{\frac12}\frac{r^{\al+\mu+1}\,e^{-\frac{(ry)^2}4}}{(\ln\frac1r)^{\nu+1}}\,\frac{dr}r.\label{aux2_L2}
\Ee
Changing variables $ry=u$, the latter integral can be estimated by
\[
y^{-(\al+\mu+1)}\,
\int_{0}^{\frac y2}\frac{u^{\al+\mu+1}\,e^{-\frac{u^2}4}}{(\ln\frac yu)^{\nu+1}}\,\frac{du}u\;
\approx\;\frac1{y^{\al+\mu+1}\,[\log y]^{1+\nu}}
\]
since the major contribution happens when $u\approx1$ (with the usual convention of
reducing a log power if $\al+\mu+1=0$).
Inserting this into \eqref{aux2_L2} (and using $y\geq2$ and $x\leq 1/2$) we obtain once again
\[
II\;\lesssim\;t^{2\nu}\,\la{x}\,\Phi(y).
\]
This concludes the proof of the case $xy\geq1$.

\sline (ii) \emph{Case $xy\leq1$:} this time $r_0(xy)=1-\frac{xy}2\geq \frac12$, so we may split
\[
B_t(x,y)\,=\;\int_0^{\frac12}\cdots\;+\;\int_{\frac12}^{r_0(xy)}\cdots\;=\;B_1\;+\,B_2\]
The first term can be handled essentially as in the previous case. Namely, if $y\leq2$
we use a similar bound to \eqref{aux4_L2}
\[
B_1\,\lesssim\,t^{2\nu}\,(xy)^{\al+\frac12}\, e^{-\frac{x^2+y^2}2}\,\int_0^{\frac1{2}}
\frac{r^{\al+\mu+1}}{(\ln\frac1r)^{\nu+1}}\,\frac{dr}r\,
\approx\;t^{2\nu}\,x^{\al+\frac12}\,e^{-\frac{x^2}2}\,\la{y}\,\lesssim\,t^{2\nu}\,\la{x}\,\Phi(y).\]
If $y\geq 2$, then $x\leq\frac12$ and $\frac{2x}y\leq\frac12$, so we may split
\[
B_1\;\leq \;\int_0^{\frac{2x}y}\cdots \;+\;\int_{\frac{2x}y}^{\frac12}\cdots\;
\]
and exactly the same computations we used in \eqref{aux3_L2}, give us the bound 
\[
B_1\lesssim\;t^{2\nu}\la{x}\Phi(y).\]
Thus we are left with the integral corresponding to $B_2$, that is the range $\frac12<r<1-\frac{xy}2$.
First of all, observe that
\[
\ln\frac1r \,\approx \,1-r,\quad r\in[1/2,1]\quad \Rightarrow \quad e^{-\frac{t^2}{2\ln\frac1r}}\leq \,e^{-\frac{ct^2}{1-r}},
\]
for a suitable $c>0$.
Next, we need once again more precise bounds for the exponentials in \eqref{aux1_L2}. We claim that, if $r\in[1/2,1]$ then 
\[
e^{-\frac{1+r^2}{1-r^2}\frac{x^2+y^2}2}\,\leq \,e^{-\ga\frac{x^2+y^2}{1-r}}\,e^{-(1+\ga)\frac{x^2+y^2}2}\,
\]
for a small constant $\ga>0$. This is easily obtained using the fact that $\frac{1+r^2}{1-r^2}\geq \frac53$ in this interval. With these exponential bounds we can control the integral $B_2$ as follows
\Bea
B_2& \lesssim& t^{2\nu}\,(xy)^{\al+\frac12}\,
e^{-(1+\ga)\frac{x^2+y^2}2}\,\int_{1/2}^{1-\frac{xy}2} \frac{e^{-\frac{ct^2+\ga(x^2+y^2)}{1-r}}}{(1-r)^{\al+\nu+1}}
\,\frac{dr}{1-r}\nonumber\\
& \lesssim & \frac{t^{2\nu}\,(xy)^{\al+\frac12}}{[t^2+x^2+y^2]^{\al+\nu+1}}\,\,
e^{-(1+\ga)\frac{x^2+y^2}2}\int_0^\infty e^{-u}\,u^{\al+\nu+1}\,\frac{du}u,\label{aux7_P2}
\Eea
after changing variables {\small $u=[ct^2+\ga(x^2+y^2)]/(1-r)$}. The last integral is a finite constant (because
$\al+\nu+1>0$), so we observe three possible cases:
\Benu
\item if $y\geq 1$, we can disregard the denominator and obtain
\[B_2\lesssim t^{2\nu}\,(xy)^{\al+\frac12}\,e^{-(1+\ga)\frac{y^2}2}\,\lesssim\,t^{2\nu}\,\la{x}\,\Phi(y),\]
since the exponential decay in $y$ is actually better than $\Phi(y)$ (and also $x\leq 1$).

\item if $y\leq1$ and $\max\{x,t\}\geq1$, we can also disregard the denominator and obtain
\[B_2\lesssim t^{2\nu}\,\la{y}\,x^{\al+\frac12}\,e^{-\frac{x^2}2}\,\lesssim\,t^{2\nu}\,\la{x}\,\la{y}.\]

\item if all $y,t,x\leq 1$, we bound the denominator in the two obvious ways to obtain
\Be
B_2\,\lesssim\,\frac{t^{2\nu}\,(xy)^{\al+\frac12}}{t^{2\nu}\,x^{2(\al+1)}}\,=\,\frac{\la{y}}{\lat{x}},\label{aux8_P2}\Ee
which is precisely the upper bound stated in \eqref{B_L2}. Observe that when $x\to0$ this piece gives the largest contribution to $B_t(x,y)$. 
\Eenu\ProofEnd

\BR\label{Rem2}
It is also possible to obtain a bound\Be
B_t(x,y)\,\lesssim\,c'(x)\,t^{2\nu}\,\Phi(y),\label{R2}
\Ee
with perhaps a worse function $c'(x)$, but without the loss produced by $\max\{1,t^{2\nu}\}$.
This loss appeared when $t,x,y\leq 1$ in \eqref{aux8_P2} above.
Looking at \eqref{aux7_P2} we may replace that bound by \[
B_2\,\lesssim\,\frac{t^{2\nu}\,(xy)^{\al+\frac12}}{x^{2(\al+\nu+1)}},\]
which implies \eqref{R2} with $c'(x)=1/\l{x}^{\al+\frac32+2\nu}$. This estimate will also be useful later.
\ER

\subsection{Upper estimates for $A_t(x,y)$: integrals over $r\leq 1/2$}
Recall that
\Be
A_t(x,y)\,\approx\,t^{2\nu}\,e^{\frac{x^2-y^2}2}\;\int_{r_0(xy)}^1
\frac{r^{\mu+\frac12}}{\sqrt{1-r}\,(\ln\frac1r)^{1+\nu}}\,
e^{-\frac{t^2}{2\ln\frac1r}}\,e^{-\frac{(x-ry)^2}{1-r^2}}\,
\frac{dr}{r}\,.
\label{Atxy}
\Ee
When $xy\geq1$ we have $r_0(xy)=\frac1{2xy}\leq \frac12$, so we can write
\[
A_t(x,y)\,=\,\int_{r_0(xy)}^{\frac12}\cdots \;+\;\int_{\frac12}^1\cdots\;=\,A1\,+\,A2.
\]
In this section we shall prove the following estimate for $A1$.

\begin{lemma}\label{L3}
If $xy\geq1$, then
\Be
A1\,\lesssim\,t^{2\nu}\,e^{\frac{x^2-y^2}2}\;\int_{\frac1{2xy}}^\frac12
\frac{r^{\mu+\frac12}}{(\ln\frac1r)^{1+\nu}}\,e^{-\frac{(x-ry)^2}{1-r^2}}\,
\frac{dr}{r}\,
\,\lesssim\,c(x)\,t^{2\nu}\,\Phi(y),\label{L3f}\Ee
where $c(x)=\,[\log(e+x)]^{\nu+1}\,(1+x)^{|\mu+\frac12|}\,\exp(x^2/2)$.
\end{lemma}
\Proof We shall distinguish cases

\sline (i) \emph{Case $y\leq 4x$}. 
In this region we essentially disregard the exponential term $e^{-\frac{(x-ry)^2}{1-r^2}}$ inside the integral, and directly estimate
\Be
A1\,\lesssim\,t^{2\nu}\,e^{\frac{x^2-y^2}2}\;\int_{\frac1{2xy}}^\frac12
\frac{r^{\mu+\frac12}}{(\ln\frac1r)^{1+\nu}}\,
\frac{dr}{r}\,.
\label{aux4_L3}\Ee
Notice however that when $y\leq x$ the exponential produces an additional gain, due to\Be
ry\leq \tfrac x2\quad\Rightarrow\quad|x-ry|\,\geq\,\tfrac{x}2\quad\Rightarrow\quad
e^{-\frac{(x-ry)^2}{1-r^2}}\,\leq\,e^{-\frac{x^2}4}.
\label{expgain}\Ee
This will play a role later in evaluating the constant $c(x)$.
We now evaluate the integral in \eqref{aux4_L3}, depending on the sign of $\mu+\frac12$.
\Benu
\item If $\mu+\frac12\geq0$ the integral is bounded by a constant, and hence
\[
A1\,\lesssim\,t^{2\nu}\,e^{\frac{x^2-y^2}2}.
\]
We shall enlarge this value to match \eqref{L3f} as follows.
Since $xy\geq1$, in this range we have $x\geq1/2$. So if $1\leq y\leq 4x$ we may use 
\[
1\,\lesssim\,\frac{(1+x)^{\mu+\frac12}[\log(e+x)]^{\nu+1}}{(1+y)^{\mu+\frac12}[\log(e+y)]^{\nu+1}}.\]
If $\frac1x\leq y\leq 1$, we use instead\[
1\,\lesssim\,\max\{x^{\al+\frac12},1\}\,\la{y},
\]
which in this region can be combined with the extra exponential in factor in \eqref{expgain}. In both cases
we obtain $A1\lesssim  t^{2\nu}c(x)\Phi(y)$, as wished.

\item If $\mu+\frac12<0$ the integral diverges near 0, but we still obtain
\[
\int_{\frac1{2xy}}^\frac12
\frac{r^{\mu+\frac12}}{(\ln\frac1r)^{1+\nu}}\,
\frac{dr}{r}\,\lesssim\,\frac{1}{(xy)^{\mu+\frac12}(\log 2xy)^{\nu+1}}.\]
Thus, using the inequality $\log(2xy)\gtrsim \max\{\log y,\log 2\}$ we arrive at 
\Beas A1 &\lesssim &t^{2\nu}\,e^{\frac{x^2-y^2}2}\,\frac{(1+x)^{|\mu+\frac12|}}
{(1+y)^{\mu+\frac12}\,[\log(e+y)]^{\nu+1}}\;\leq\;t^{2\nu}c(x)\Phi(y),\quad \mbox{if }\;1\leq y\leq 4x\\
 &\lesssim &t^{2\nu}\,\frac{e^{\frac{x^2}4}\,y^{\al+\frac12}}
{x^{\mu+\frac12}\,y^{\mu+\al+1}}\,\leq\,t^{2\nu}\,e^{\frac{x^2}4}x^{\al+\frac12}\,\la{y}\;\leq\;t^{2\nu}c(x)\Phi(y),\quad \mbox{if }\,\frac1x\leq y\leq1,
\Eeas
using in the last case the additional exponential gain in \eqref{expgain}. 
\Eenu

\sline (ii) \emph{Case $y\geq 4x$}. This is the same as $\frac{2 x}y\leq \frac12$, and remember from \eqref{aux6_L2}
that when $r\in[\frac {2x}y,\frac12]$ a better bound for the exponential is available, namely
\Be
e^{-\frac{(x-ry)^2}{1-r^2}}\,\leq\,e^{-\frac{r^2y^2}4}.\label{aux2_L3}\Ee
Thus we may consider two subcases, depending on whether $\frac{2x}y$ is above or below $r_0(xy)$.

\sline $\bullet$ \emph{Subcase $\frac{2x}y\leq r_0(xy)=\frac1{2xy}\leq \frac12$.} Using \eqref{aux2_L3} we obtain
\Beas
A1& \lesssim & t^{2\nu}\,e^{\frac{x^2-y^2}2}\;\int_{\frac1{2xy}}^\frac12
\frac{r^{\mu+\frac12}\,e^{-\frac{(ry)^2}4}}{(\ln\frac1r)^{1+\nu}}\,
\frac{dr}{r}\,\\
\mbox{{\footnotesize ($ry=u$)}}& = & \frac{t^{2\nu}\,e^{\frac{x^2-y^2}2}}{y^{\mu+\frac12}}
\;\int_{\frac1{2x}}^\infty
\frac{u^{\mu+\frac12}\,e^{-\frac{u^2}4}}{(\ln\frac yu)^{1+\nu}}\,
\frac{du}{u}\,\Eeas
Observe that $x\leq \frac12$ (and $y\geq2$), so the latter integral reaches its major contribution at $u=\frac1{2x}$
\[
\int_{\frac1{2x}}^\infty
\frac{u^{\mu+\frac12}\,e^{-\frac{u^2}4}}{(\ln\frac yu)^{1+\nu}}\,
\frac{du}{u}\,\lesssim \frac{(1/x)^{\mu-\frac12}e^{-c/x^2}}{(\log 2xy)^{1+\nu}}\;\lesssim\;\frac1{(\log y)^{1+\nu}},\]
using in the last step the elementary bound of logarithms
\[
\log(2xy)\,\gtrsim\,\frac{\log(y+e)}{\log(\frac1x+e)}, \quad \mbox{if $y\geq\max\{1,1/x\}$}
\]
(see e.g. \cite[Lemma 5.1]{gar}). Thus we conclude that\[
A1\,\lesssim \,t^{2\nu}\,e^{\frac{x^2}2}\,\Phi(y).\]

\sline $\bullet$ \emph{Subcase $r_0(xy)<\frac{2x}y\leq \frac12$.} Here we split
\[
A1\,= \,\int^\frac12_{\frac{2x}y}\cdots \;+\;\int^{\frac{2x}y}_{r_0(xy)}\cdots\;=I\,+\,II.
\]
The first term is similar to the previous subcase, except that now $x>\frac12$ (and $y\geq 4x\geq 2$)
\Beas
I & \lesssim & 
\frac{t^{2\nu}\,e^{\frac{x^2-y^2}2}}{y^{\mu+\frac12}}\;\int_{\frac x2}^\infty
\frac{u^{\mu+\frac12}\,e^{-\frac{u^2}4}}{(\ln\frac yu)^{1+\nu}}\,
\frac{du}{u}\,\Eeas
and the last integral is bounded by a constant times
\[
\frac{x^{\mu-\frac12}\,e^{-cx^2}}{(\ln\frac {2y}x)^{1+\nu}}\,
\lesssim\,x^{\mu-\frac12}\,e^{-cx^2}\,\frac{[\log(e+x)]^{1+\nu}}{[\log(e+y)]^{1+\nu}}\,
\lesssim\,\frac{1}{[\log(e+y)]^{1+\nu}}.\]
Finally, we consider $II$. Here there is no exponential gain, and similarly to \eqref{aux4_L3}
we have
\[
II\,\lesssim\,t^{2\nu}\,e^{\frac{x^2-y^2}2}\;\int_{\frac1{2xy}}^\frac{2x}y
\frac{r^{\mu+\frac12}}{(\ln\frac1r)^{1+\nu}}\,
\frac{dr}{r}\,=\,\frac{t^{2\nu}\,e^{\frac{x^2-y^2}2}}{y^{\mu+\frac12}}\;\int_{\frac1{2x}}^{2x}
\frac{u^{\mu+\frac12}}{(\ln\frac  yu)^{1+\nu}}\,
\frac{du}{u}\,.\]
Now, the last integral can easily be analyzed (depending on the sign of $\mu+\frac12$) to obtain \[
\int_{\frac1{2x}}^{2x}
\frac{u^{\mu+\frac12}}{(\ln\frac  yu)^{1+\nu}}\,
\frac{du}{u}\,\lesssim\,\frac{x^{|\mu+\frac12|}[\log(x+e)]^{1+\nu}}{[\log(y+e)]^{1+\nu}}.
\]
Thus, overall we conclude that in this subcase
\[
A1\;\lesssim\;I\,+\,II\;\lesssim\,t^{2\nu}\,(1+x)^{|\mu+\frac12|}[\log(x+e)]^{1+\nu}\,e^{\frac{x^2}2}\,
\Phi(y).
\]
\ProofEnd

\subsection{Upper estimates for $A_t(x,y)$ when $y\leq x/2$ or $y\geq Mx$}

In view of the previous subsection, it only remains to estimate
\Bea
A2& \approx & t^{2\nu}\,e^{\frac{x^2-y^2}2}\;\int_{\max\{r_0(xy),\frac12\}}^1
\frac{e^{-\frac{t^2}{2\ln\frac1r}}}{\sqrt{1-r}\,(\ln\frac1r)^{1+\nu}}\,
\,e^{-\frac{(x-ry)^2}{1-r^2}}\,{dr}\,\nonumber\\
& \lesssim & t^{2\nu}\,e^{\frac{x^2-y^2}2}\;\int_{\max\{1-\frac{xy}2,\frac12\}}^1
\frac{e^{-\frac{ct^2}{1-r}}}{(1-r)^{\nu+\frac32}}\,
\,e^{-\frac{(x-ry)^2}{1-r^2}}\,{dr}\,\label{aux1_L4}\Eea
noticing that $\ln\frac1r\approx 1-r$ when $r\in[\frac12,1]$. In this region, however, it is more convenient to use
the write-up for the heat kernel in \eqref{hs}, in terms of the parameter $s$. This gives a more reasonable 
expression for the exponentials, namely\[
e^{\frac{x^2-y^2}2}\,e^{-\frac{(x-ry)^2}{1-r^2}}\,=
\,e^{-\frac14[\frac{(x-y)^2}{s}+s(x+y)^2]}.
\]
Since the parameters $r$ and $s$ are related by $s=\frac{1-r}{1+r}$ (or $r=\frac{1-s}{1+s}$), 
either from \eqref{hs} or directly from \eqref{aux1_L4}, we obtain that
\Be
A2\;\lesssim\;t^{2\nu} \int_0^{\min\{\frac13,\frac{xy}3\}}\frac{e^{-\frac{ct^2}{s}}\,e^{-\frac14[\frac{(x-y)^2}{s}+s(x+y)^2]}}{s^{\nu+\frac32}}\,ds,\label{As}\Ee
after perhaps slightly enlarging the range of integration.
Our first result shows that when $y$ is far from $x$ this can also be controlled by the function $\Phi(y)$.

\begin{lemma}\label{L4}
There exists $M>1$ such that, if $\;y\leq \frac x2\;$ or $\;y\geq Mx\;$, then
\[
A2\;\lesssim\;t^{2\nu} \int_0^{\min\{\frac13,\frac{xy}3\}}\frac{e^{-\frac{ct^2}{s}}\,e^{-\frac14[\frac{(x-y)^2}{s}+s(x+y)^2]}}{s^{\nu+\frac32}}\,ds\;\lesssim\; c(x)\,\max\{t^{2\nu},1\}\,\Phi(y),
\]
with $c(x)=1/\lat{x}$.
\end{lemma}
\Proof  We claim that, in the assumed range of $x$ and $y$, there is  some $\ga>0$ such that
\Be
A2\;\lesssim\;t^{2\nu}\,e^{-(\frac12+\ga)y^2}\; \int_0^{\min\{\frac13,\frac{xy}3\}}{e^{-\ga\frac{t^2+(x-y)^2}{s}}}
{s^{-(\nu+\frac32)}}\,ds
\label{aux2_L4}\Ee
This is just a bound of the exponentials. Indeed, if we distinguish the two cases

\Benu\item \emph{case $y\geq Mx$}: this implies $|y-x|\geq (1-\frac{1}M)y$, so for any $\eta<1$ we have 
\[
e^{-\frac{ct^2}{s}}\,e^{-\frac14[\frac{(x-y)^2}{s}+s(x+y)^2]}\,\leq\,
e^{-\frac{ct^2+\frac\eta4(x-y)^2}{s}}\,e^{-\frac{1-\eta}4(\frac{M-1}M)^2(\frac{1}{s}+s)y^2}\,
\]
which implies the required assertion using that $\frac1{s}+s\geq\frac{10}{3}$, when $s\in(0,\frac13)$,
 and choosing $M$ sufficiently large and $\eta$ sufficiently
small.

\item \emph{case $y\leq x/2$}: this time $|x-y|\geq\frac x2\geq y$, so we have 
\[
e^{-\frac{ct^2}{s}}\,e^{-\frac14[\frac{(x-y)^2}{s}+s(x+y)^2]}\,\leq\,
e^{-\frac{ct^2+\frac\eta4(x-y)^2}{s}}\,e^{-\frac{1-\eta}4(\frac{1}{s}+s)y^2}\,
\]
which again implies the assertion using $\frac1{s}+s\geq\frac{10}{3}$
 and choosing $\eta$ sufficiently
small.
\Eenu
Thus \eqref{aux2_L4} is proven, and we may change variables $[t^2+(x-y)^2]/s=u$ to obtain
\Bea
A2 & \lesssim & \frac{t^{2\nu}\,e^{-(\frac12+\ga)y^2}}{[t^{2}+(x-y)^2]^{\nu+\frac12}}
\; \int^\infty_{3\ga[t^2+(x-y)^2]\max\{1,\frac1{xy}\}}
{u^{\nu+\frac12}}\,{e^{- u}}\,\frac{du}u\label{aux3a_L4} \\
& \lesssim & \frac{t^{2\nu}\,e^{-(\frac12+\ga)y^2}}{(t+x+y)^{2\nu+1}}
\; \int^\infty_{\ga'[x^2+y^2]\max\{1,\frac1{xy}\}}
{u^{\nu+\frac12}}\,{e^{- u}}\,\frac{du}u,\label{aux3_L4}\Eea
since in the selected range of $x,y$ we have $|x-y|\gtrsim x+y$. To finish the proof we
must distinguish some cases.

\sline \emph{Case $y\geq1$:} then bounding the denominator and the integral in \eqref{aux3_L4} by a constant
we immediately see that
\[
A2\;\lesssim\;t^{2\nu}\,e^{-(\frac12+\ga)y^2}\,\lesssim\,t^{2\nu}\,\Phi(y)\]
since the exponential has a better decay.

\sline \emph{Case $y\leq1$ and $y\geq Mx$:} we again bound the integral by a constant and 
estimate the fraction in \eqref{aux3_L4} as follows
\Be
A2\;\lesssim\;\frac{t^{2\nu}}{t^{2\nu}y}=\frac{\la{y}}{y^{\al+\frac32}}\,\leq\,c_M\, \frac{\la{y}}{\lat{x}}.
\label{aux8_L4}\Ee

\sline \emph{Case $y\leq1$ and $y\leq \frac x2$:} this is a relevant case, since the integral in \eqref{aux3_L4} plays actually a role. To  compute the integral we must distinguish the two subcases

\Benu
\item If $xy\leq1$, then since also $\frac xy\geq 2$,
\Bea
A2 & \lesssim &\frac{t^{2\nu}}{t^{2\nu}x}\;\int^\infty_{\ga'\frac x{y}}
{u^{\nu-\frac12}}\,{e^{- u}}\,{du}\,\approx\,\frac1x\,\Big(\frac xy\Big)^{\nu-\frac12}e^{-c\frac xy}
\label{aux9_L4}\\
&\lesssim&\frac 1x\,\Big(\frac yx\Big)^{\al+\frac12}\;\leq\;  \frac{\la{y}}{\lat{x}}.\nonumber
\Eea

\item If $xy\geq1$, then we have  $x\geq \frac 1y\geq 1$ and
\[
A2 \; \lesssim \;\frac{t^{2\nu}}{x^{1+2\nu}}\;\int^\infty_{\ga'\,x^2}
{u^{\nu-\frac12}}\,{e^{- u}}\,{du}\;\lesssim\;t^{2\nu}\,x^{-2}\,e^{-c\,x^2}.\]
Now, since $\frac1x\leq y \leq 1$ we can insert the estimate\[
1\,\lesssim\, \la{y}\,\max\{x^{\al+\frac12},1\},
\]
to obtain $A2\lesssim t^{2\nu}\la{y}$.
\Eenu
\ProofEnd
\BR
\label{Rem3}
As mentioned earlier in Remark \ref{Rem2}, here it is also possible to obtain a bound\Be
A2\,\lesssim\,c'(x)\,t^{2\nu}\,\Phi(y),\label{R3}
\Ee
with $c'(x)=1/\l{x}^{\al+\frac32+2\nu}$.
The loss produced by $\max\{1,t^{2\nu}\}$ can be corrected in 
\eqref{aux8_L4} and \eqref{aux9_L4} by replacing the factor $t^{2\nu}$ in the denominator
by $x^{2\nu}$, as one readily notices from \eqref{aux3_L4}. 
As mentioned before, this estimate will play a role later.
\ER

\subsection{Upper estimates for $A_t(x,y)$ in the local part $\frac x2<y<Mx$}

As in the previous subsection, our starting point is the formula \eqref{As},
which we must estimate in the  local region $\frac x2<y<Mx$. A sufficient bound for us is stated
in the next lemma.

\begin{lemma}\label{L5}
If $\;\frac x2<y<Mx\;$, then
\Be
t^{2\nu} \int_0^{\min\{\frac13,\frac{xy}3\}}\frac{e^{-\frac{ct^2}{s}-\frac14[\frac{(x-y)^2}{s}+s(x+y)^2]}}{s^{\nu+\frac32}}\,ds\;\lesssim\;C(x)\;\frac{ t^{2\nu}\,e^{-\frac{y^2}2}}{\big(t+|x-y|\big)^{1+2\nu}},\label{L5a}
\Ee
where  $C(x) \,=\,(1+x)^{2\nu}e^{\frac{x^2}2}$.
\end{lemma}
\Proof
We shall crudely enlarge the integral in \eqref{L5a} to the range $\int_0^{1/3}$. 
This last integral was already estimated in \cite{GHSTV} and \cite{gar}, by a similar procedure to the one used in the last subsection. 
More precisely, from the estimates in \cite[Lemma 3.2]{gar}, formula (3.16), it follows that
\[
 t^{2\nu} \int_0^{\frac12}\frac{e^{-\frac{ct^2}{s}-\frac14[\frac{(x-y)^2}{s}+s(x+y)^2]}}{s^{\nu+\frac32}}\,ds\;\lesssim\;\frac{t^{2\nu}\,(1+x)^{2\nu}\,e^{\frac{x^2-y^2}2}}{\big(t+|x-y|\big)^{1+2\nu}}
,\]
which agrees with \eqref{L5a}. 
\ProofEnd

\subsection{Proof of Propositions \ref{P3.1} and \ref{P3.2}}

Proposition \ref{P3.2} follows by putting together Lemmas \ref{L2}, \ref{L3}, \ref{L4} and \ref{L5}.
Concerning Proposition \ref{P3.1}, the lower bound was shown in Lemma \ref{L1}, while
the upper bound also follows from Lemmas \ref{L2}, \ref{L3} and \ref{L4}, at least when $y<\frac x2$ or $y>Mx$.
This actually implies the asserted result for all $x$ and $y$, since when $y$ belongs to the compact set $[\frac x2,Mx]$, the continuous function $y\to P_t(x,y)/\Phi(y)$ is bounded above by a constant $c_2(t,x)$.
\ProofEnd

\subsection{Proof of Corollary \ref{C3.3}}

By Proposition \ref{P3.2}, observe that\Be
P_tf(x)\,\lesssim\,\frac{C_1(x)}{t}
\,\int_{\SRp}\frac{g(y)\,dy}{(1+\frac{|x-y|}t)^{1+2\nu}}\;+\;
C_2(x)(1\vee t_0)^{2\nu}\,\big\|f\big\|_{L^1(\Phi)},\label{aux1_C3}\Ee
where $g(y)=f(y)e^{-\frac{y^2}2}\chi_{\{\frac x2<y<Mx\}}$.
The first term is then controlled by a maximal function by a standard slicing argument.
\ProofEnd

\BR
We wrote in \eqref{Pcontrol} a different version of \eqref{P**} with $\Ml(f\Phi)$ 
in place of $\Ml(fe^{-\frac{y^2}2})$. Since $x\approx y$, 
\[
\Ml(f\Phi)(x)\,\approx\, \tfrac{\la{x}}{[\log(e+x)]^{1+\nu}\,(1+x)^{\mu+\frac12}}\,\Ml(fe^{-\frac{y^2}2})(x),
\] so they are actually equivalent modulo $x$-constants. 
 The write-up in \eqref{Pcontrol}
has the advantage of remaining valid for other Laguerre systems; see $\S\ref{transf}$ below.
\ER

\section{Proofs}\label{proofs}

As indicated in $\S1$ we present the proof of Theorems \ref{th1SL} and \ref{th2SL}
for the differential operator $L$ in \eqref{Lphi} and the function $\Phi$ in \eqref{phiL}.
We postpone to $\S\ref{transf}$ the proof of the results for the other systems mentioned in Table 1.

\subsection{Proof of Theorem \ref{th1SL}}
First of all, it is an immediate consequence of Proposition \ref{P3.1} that $P_t|f|(x)<\infty$  for some (or all) 
$t,x>0$ if and only if $f\in L^1(\Phi)$. This justifies that $f\in L^1(\Phi)$
is the right setting for this problem. Notice also that taking derivatives of the kernel $P_t(x,y)$
in \eqref{Ptxy} with respect to $t$ does not worsen its decay in $y$, so $P_tf(x)$ automatically becomes infinitely differentiable
in the $t$-variable when $f\in L^1(\Phi)$. We can also take as many derivatives as wished with respect to $x$,
since the kernel satisfies the pde\footnote{For a justification that the subordinated integral in \eqref{poisf} satifies the pde \eqref{pde}, see e.g. \cite[$\S2$]{gar}.} 
\[\Big[\partial_{xx} 
- x^2 - \tfrac{\al^2-\tfrac14}{x^2}\,-\,2\mu\,+\,\partial_{tt} \,+\,\tfrac{1-2\nu}t\,\partial_t\Big]P_t(x,y)=0,\]
so $x$-derivatives are transformed into $t$-derivatives and 
do not worsen the decay of $P_t(x,y)$ in the $y$-variable.  
We have thus completed the proof of paragraph (i) and the last statement in Theorem \ref{th1SL}. We shall now prove a stronger result than (ii), namely that for $f\in L^1(\Phi)$
\Be
\lim_{t\to0^+}P_tf(x)=f(x), \quad\forall\;x\in\cL_f
\label{ptwise}\Ee
where $\cL_f$ denotes the set of Lebesgue points of $f$. When $f(x)=0$ this is
easily obtained from the kernel estimates in Proposition \ref{P3.2}. Indeed,
\[
P_tf(x) \,\lesssim\,C_1(x)\int_{\frac x2}^{Mx}\frac{t^{2\nu}|f(y)|\,dy}{(t+|x-y|)^{1+2\nu}}\,+\,
C_2'(x)\,t^{2\nu}\,\int_{\SR_+}|f|\Phi,
\]
where in the second term we have replaced $(t\vee 1)^{2\nu}$ by $t^{2\nu}$ in view of 
 Remarks \ref{R2} and \ref{R3}. Thus, this second term vanishes as $t\to0$ (actually for all $x\in\SR_+$).
Concerning the first term, it is given by convolution of $|f(y)|\chi_{\{\frac x2<y<Mx\}}\in L^1_c(\SR_+)$
with a radially decreasing approximate identity, so from well-known results (see e.g. \cite[p. 112]{SteSha}), it must vanish as $t\to0$ at every Lebesgue point $x$ of $f$ with $f(x)=0$.

It remains to prove \eqref{ptwise} when $f(x)$ is not necessarily 0. To show this, we first notice that
the first eigenfunction $\phi=\pz$ (with eigenvalue $\lam=2(\mu+\al+1)$) satisfies \[P_t\phi=F_t(\lam)\phi,\quad \mbox{with}\quad  \lim_{t\to0}F_t(\lam)=1.\]
Indeed, setting $u={t^2}/{(4v)}$ in \eqref{poisf}, gives
$$F_t(\lam)=\tfrac{(t/2)^{2\nu}}{\Ga(\nu)}\,\int_0^\infty 
e^{-\frac{t^2}{4u}-\lam u}\,\frac{du}{u^{1+\nu}}=\tfrac{1}{\Ga(\nu)}\int_0^\infty e^{-v-\frac{t^2\lam}{4v}}\,v^{\nu-1}\,dv\longrightarrow1,\quad \mbox{as }t\to0.$$
Therefore, we can write
\Be
P_tf(x)-f(x)=\,P_tf(x)-F_t(\lam)f(x) \;+\;f(x)[F_t(\lam)-1],
\label{aux_Ft}\Ee
with the last term vanishing as $t\to0$. Since $\phi>0$, the first term can be rewritten as\[
P_tf(x)-\tfrac{f(x)}{\phi(x)}\,P_t\phi(x)\,=\, P_t\left(f-\tfrac{f(x)}{\phi(x)}\phi\right)(x).
\]
Setting $g=f-\tfrac{f(x)}{\phi(x)}\phi$, it is easily seen that $g\in L^1(\Phi)$, $g(x)=0$ and
$x$ is a Lebesgue point of $g$. This last assertion follows from \[
\mint_{I_r(x)}|g(y)|dy \leq \mint_{I_r(x)}|f(y)-f(x)|dy +\tfrac{|f(x)|}{\phi(x)}\mint_{I_r(x)}|\phi(y)-\phi(x)|dy,
\]
which vanishes as $r\to0$. Thus we can apply our earlier case to $g$ and conclude that $\lim_{t\to0}P_tg(x)=0$. So the left hand side of \eqref{aux_Ft} goes to 0 as $t\to0$, establishing \eqref{ptwise}
and completing the proof of
Theorem \ref{th1SL}.
\ProofEnd

\BR
A close look at the last part of the proof shows that, when $f\in C([a,b])$ with $[a,b]\Subset\SR_+$,
then the convergence of $P_tf(x)\to f(x)$ is uniform in $x\in[a,b]$. 
\ER

\subsection{Proof of Theorem \ref{th2SL}}\label{proofth2}
We have to show that $P^*_{t_0}$ maps $L^p(w)\to L^p(v)$, for some weight $v(x)>0$,
under the assumption that
\[
\|w\|_{D_p(\Phi)}:=\,\Big[\int_{\SRp}\wp(x)\Phi(x)^{p'}\,dx\Big]^{1/p'}\,<\,\infty,\]
with $\Phi$ defined as in \eqref{phiL}. We shall use the bound for $P^*_{t_0}$ in \eqref{P**},
namely \Bea
P^*_{t_0}f(x) & \lesssim &  C_1(x)\,\Ml_M \big(fe^{-\frac{y^2}2}\big)(x) \,+\,
C_2(x)\,\int_{\SRp}|f(y)|\,\Phi(y)\,dy\nonumber\\ & = & \,
I(x)\,+\,II(x),\label{IyII}\Eea
for a suitable $M>1$, and $C_1(x), C_2(x)$ given explicitly in Proposition \ref{P3.2}.
We first treat the last term, which  by H\"older's inequality is bounded by
\[
II(x)\,\leq\, C_2(x)\,\|f\|_{L^p(w)}\,\|w\|_{D_p(\Phi)}\;.
\]
Thus, it suffices to choose a weight $v$ such that $C_2(x)= [\log(e+x)]^{1+\nu}\,{(1+x)^{|\mu+\frac12|}\,e^{\frac{x^2}2}}/{\lat{x}}$ belongs to $L^p(v)$ to conclude that \Be
\|II\|_{L^p(v)} \leq\,\|C_2\|_{L^p(v)}\,\|w\|_{D_p(\Phi)}\, \|f\|_{L^p(w)}\,.
\label{IIn}\Ee 
For instance we may take any $v(x)\leq v_2(x)$ with
\Be
v_2(x):=\,\frac{\l{x}^{(\al+\frac32)p-1}}{[\log(e/\l{x})]^2}\,\frac{e^{-\frac
p2x^2}}{(1+x)^{N}}\,
\label{v2p}\Ee
for any $N>1+p\,|\mu+\frac12|$. We remark that $v_2\in D_q(\Phi)$ for all $q>p$.
Indeed, the local condition was already established in \eqref{Dqv1}.
For the global condition notice that\[
\int_1^\infty v_2(x)^{-\frac{q'}q}\,\Phi(x)^{q'}\,dx \,\lesssim\,
\int_1^\infty e^{-(1-\frac pq)\frac
{q'}2x^2}\,(1+x)^{\frac{Nq'}q}\,dx\,<\,\infty.
\]
We now consider the term $I(x)$ in \eqref{IyII}. We define a new weight $W(x)=w(x)e^{\frac p2x^2}$, and observe that
\Be
w\in D_p(\Phi)\quad \Longrightarrow\quad W\in D^0_p(\al+\tfrac12)\cap D^{\rm exp}_p(a),\quad \forall\;a>0,
\label{wWDp}\Ee
for the weight classes defined just before Proposition \ref{P5.2}.
Indeed, the local estimate follows from $\Phi(x)\approx\la{x}$ when $x\in(0,1)$, and the global estimate is a consequence of
\Be
\big\|W\big\|_{D^{\rm exp}_p(a)} \, = \, \int_1^\infty\wpp(x)\,e^{-\frac{p'}2x^2}\,
e^{-ap'x^2}\,dx\,\leq\, C^{p'}_{a}\,\int_1^\infty\wpp(x)\Phi(x)^{p'}\,dx\,,
\label{WDexp}\Ee
with $C_a=\max_{x\geq1}|\log(e+x)|^{1+\nu}|1+x|^{\mu+\frac12}e^{-ax^2}<\infty$. 

\

We shall now set\Be
v_{1,\e}(x)=\frac{e^{-\frac{px^2}2}}{(1+x)^{2p\nu}}\,
\mathscr{V}(x)\rho_{\e}\Big({\mathscr{V}(x)}\Big),\quad\mbox{where }
\mathscr{V}(x)=\Big[\Ml_M\big(\wpp e^{-\frac{p'x^2}2}\big)(x)\Big]^{-\frac p{p'}}
\label{v1p}\Ee
(or $v_{1,\e}(x)=(1+x)^{-2p\nu}e^{-\frac{px^2}2}V_{\e}(x)$  in the notation of \eqref{Ve}).
Given $f\in L^p(w)$, we denote $\tf(y)=f(y)e^{-\frac{y^2}2}$ which is a function in $L^p(W)$.
Then, using the two-weight inequality for $\Ml_M$ in Theorem \ref{ThVe}, and the expression 
for $C_1(x)=(1+x)^{2\nu}e^{\frac{x^2}2}$, we see that, for any $v\leq v_{1,\e}$, the term $I(x)$ in \eqref{IyII} is controlled by
\Bea
\|I(x)\|^p_{L^p(v)} &  \leq & \int_{\SRp}\frac{C_1(x)^p\,e^{-\frac p2x^2}}{(1+|x|)^{p2\nu}}\,
\big|\Ml_M\tf(x)\big|^p\,V_{\e}(x)\,dx\nonumber\\
& \lesssim & \|\tf\|_{L^p(W)}^p\,=\,\|f\|_{L^p(w)}^p.
\label{In}\Eea
So, combining \eqref{IyII}, \eqref{IIn} and \eqref{In} we have shown that 
$\|P^*_{t_0} f\|_{L^p(v)}\lesssim\|f\|_{L^p(w)}$,
provided \Be v(x)= \min\{v_{1,\e}(x),v_2(x)\},\label{vdef_pois}\Ee
with $v_{1,\e}(x)$ and $v_2(x)$ defined in \eqref{v1p} and \eqref{v2p}.

We only have to verify that, if $q>p$ then we can choose $\e$ sufficiently small
so that $v_{1,\e}\in D_q(\Phi)$ (which implies $v\in D_q(\Phi)$). 
This actually follows from \eqref{wWDp} and Proposition \ref{P5.2}. Indeed, on the one hand,
since $W\in D^0_p(\al+\frac12)$ \Be
\int_0^1 v_{1,\e}(x)^{-\frac{q'}q}\,\Phi(x)^{q'}\,dx\,\lesssim\,
\int_0^1 V_{\e}(x)^{-\frac{q'}q}\,\l{x}^{(\al+\frac12)q'}\,dx,
\label{aux9a}\Ee
which is finite by (i) in the proposition (choosing $\e$ sufficiently small).
On the other hand, \Be
\int_1^\infty v_{1,\e}(x)^{-\frac{q'}q}\,\Phi(x)^{q'}\,dx\,\lesssim\,
\int_1^\infty V_{\e}(x)^{-\frac{q'}q}\,e^{-q'(1-\frac pq)\frac{x^2}2}\,
(1+x)^{2\nu pq'/q}\,dx\,,\label{aux9}\Ee
and since $W\in D^{\rm exp}_p(a)$ for all $a>0$, we can apply
part (ii) of Proposition \ref{P5.2}  
(for a sufficiently small $\e$) to conclude that this is also a finite quantity.
\ProofEnd

\BR \label{vPhi}{\bf Alternative expression for the second weight.} 
A slight modification of the above construction suggests to define a new weight by
\Be
v^{\Phi,w}_\e(x)\,:=\,\min\left\{\Phi(x)^p\,\Big[\Ml\big(w^{-\frac{p'}{p}}\Phi^{p'}\big)(x)\Big]^{-\frac{p}{p'}}\,\Upsilon_\e(x),\;\frac{\l{x}^{p-1}}{[\log(e/\l{x})]^2}\,\frac{\Phi(x)^p}{(1+x)^{N_0}}\right\}
\label{tvdef}\Ee
with
\[
\Upsilon_\e(x)=\, \frac{\l{x}^{\e N_1}}{(1+x)^{N_2}}\,
\rho_{\e}\Big(\Big[\Ml\big(w^{-\frac{p'}{p}}\Phi^{p'}\big)(x)\Big]^{-\frac{p}{p'}}\Big).
\]
If $N_0,N_1,N_2$ are sufficiently large, then similar arguments as above lead to the boundedness
of $P^*_{t_0}:L^p(w)\to L^p(v^{\Phi,w}_\e)$ for all $\e>0$, and give also the property that $v^{\Phi,w}_\e\in D_q(\Phi)$ if $\e$ is sufficiently small. We omit the details. 
The expression in \eqref{tvdef} has the advantage of remaining valid for the other Laguerre systems in Table \ref{table1}
(with the corresponding $\Phi$ functions).
\ER

\section{Transference to other Laguerre type systems}\label{transf}

In this section we show how to transfer the results already proved for the system $\{\pn\}$ 
and the operator $L$ to the other Laguerre systems and operators in Table \ref{table1}. 
The procedure is completely general, as one can infer already from the first two cases.

\subsection{Results for the system $\psi^\al_n$}
\label{transf_psi}

The starting point is the identity defining $\psi^\al_n$, namely\Be
\psi^\al_n(y)\,=\,a(y)\,\pn(y),\quad\mbox{with } a(y)=y^{-\al-\frac12}.\label{apsi}\Ee
Clearly, $\pn$ is an eigenvector of $L$ if and only if $\psn$ is an eigenvector of the operator 
\[f\longmapsto\La f(x)= a(x)L[a^{-1}f\;](x)\]
(with the same eigenvalue $\lam_n=4n+2(\al+1+\mu)$). 
An elementary computation shows that the differential operator $\La$ obtained in this fashion is exactly the one listed in Table \ref{table1}. Remark also that $\{\psn\}$ becomes an orthonormal basis in $L^2$ with the measure $a^{-2}(y)dy=y^{2\al+1}dy$.

The identity in \eqref{apsi} leads to a pointwise relation of the corresponding heat kernels
\[
e^{-t\La}(x,y)\,=\,\sum_{n=0}^\infty e^{-\lam_nt}\psn(x)\psn(y)\,=\,
a(x)\,a(y)\,e^{-tL}(x,y),
\]
and by the subordination formula, also of the corresponding Poisson kernels
\[
P_t^{\La}(x,y)\,=\,a(x)\,a(y)\,P_t^{L}(x,y).
\]
In particular,
\Be
P_t^\La f(x)\,=\,\int_{\SRp}P_t^\La(x,y)\,f(y)\,a^{-2}(y)dy\,=\,a(x)\,P^L_t[a^{-1}f](x).\label{PLaL}\Ee
From this relation it is clear that Theorem \ref{th1SL} becomes true for the operator $\La$ with
\[
\Phi^\La(y)=a(y)^{-1}\Phi^L(y),
\]
as listed in Table \ref{table1}. From \eqref{PLaL} it also follows that $P^{*,\La}_{t_0}$ maps $L^p(w)\to L^p(v)$
if and only if $P^{*,L}_{t_0}$ maps $L^p(a^pw)\to L^p(a^pv)$, and hence the necessary and sufficient condition becomes
\[
a^pw\in D_p(\Phi^L)\Longleftrightarrow\|a^{-1}w^{-\frac1p}\Phi^L\|_{p'}<\infty\Longleftrightarrow w\in D_p(\Phi^\La),
\]
as was claimed in Theorem \ref{th2SL}. For the assertions about the weight $v$ one may argue directly as follows.
Observe from \eqref{PLaL} and Corollary \ref{C3.3} that we can write
\Beas P^{*,\La}_{t_0}f(x)
& \lesssim & 
C_1(x)\,a(x)\Ml_M \big(fa^{-1}e^{-\frac{y^2}2}\big)(x)\,+\,C_2(x)a(x)\,\int_{\SRp}
|f|\Phi^\La,\Eeas
with a cancellation in the first term due to $a(x)a(y)^{-1}\approx1$ when $\frac x2<y<Mx$.
At this point we can apply the same arguments as in $\S\ref{proofth2}$. Namely, we construct 
$v=\min\{v_{1,\e},v_2\}$ with the same choice of $v_{1,\e}$, and with $v_2$ in \eqref{v2p} now replaced by
\[
v_2(x):=\,\frac{\l{x}^{(2\al+2)p-1}}{[\log(e/\l{x})]^2}\,\frac{e^{-\frac
p2x^2}}{(1+x)^{N}}\,.
\]
The same proof will give that, for any $q>p$, there is
a sufficiently small $\e$ so that $v\in D_q(\Phi^\La)$ (the only difference being that, locally, this condition now becomes $v\in D^0_q(2\al+1)$). We remark that this part will work as well with the choice\[
v(x)\,=\,a^{-p}(x)\,v^{\Phi^L, a^pw}(x)\,=\,v^{\Phi^\La,w}(x),
\]
as defined in \eqref{tvdef}.

\subsection{Results for the system $\fLn$}\label{subsecfLn}
Consider the following isometry of $L^2(\SRp,dy)$
\[
f\longmapsto Af(x)\,=\,\sqrt{2x}\,f(x^2).
\]
The systems $\fLn$ and $\pn$ are related by $\pn=A\fLn$, or equivalently\Be
\fLn(y)\,=\,[A^{-1}\pn](y)\,=\,(4y)^{-\frac14}\,\pn(\sqrt y).\label{AfLn}\Ee
In particular, $\fLn$ is an eigenvector  of the operator 
\[\fL \,= \frac14\, A^{-1}\circ L\circ A,\]
this time with eigenvalue $\lam_n/4=n+(\al+1+\mu)/2$. 
The factor $\frac14$ has been added so that $\fL$ coincides with the operator listed in Table \ref{table1}. 

The heat kernels are now related by 
\[
e^{-t\fL}(x,y)\,=\,\sum_{n=0}^\infty e^{-t\lam_n/4}\fLn(x)\fLn(y)\,=\,
\frac1{2(xy)^{\frac14}}\,e^{-\frac t4 L}(\sqrt x,\sqrt y),
\]
and therefore, a substitution in the subordinated integral in \eqref{Ptxy} gives
\[
P_t^{\fL}(x,y)\,=\,(16xy)^{-\frac14}\,P_{t/2}^{L}(\sqrt x,\sqrt y).
\]
Thus, we obtain the formula
\Be
P_t^\fL f(x)\,=\,\int_{\SRp}P_t^\fL(x,y)\,f(y)\,dy\,=\,(4x)^{-\frac14}\,P^L_{t/2}[Af](\sqrt x).\label{PfL}\Ee
From this relation one easily deduces Theorem \ref{th1SL} for the operator $\fL$, provided that
\Be
\Phi^\fL(y)=[A^{-1}\Phi^L](y),
\label{PhifL}\Ee
which is the function asserted in Table \ref{table1} (modulo constants).

To establish the second theorem, first observe from \eqref{PfL} and Corollary \ref{C3.3} that 
\Beas P^{*,\fL}_{t_0}f(x)
& \lesssim & 
\frac{C_1(\sqrt x)}{x^{1/4}}\,\Ml_M \big(\sqrt yf(y^2)e^{-\frac{y^2}2}\big)(\sqrt x)\,+\,\frac{C_2(\sqrt x)}{x^{1/4}}\,\int_{\SRp}
\sqrt{2y}|f(y^2)|\Phi^L(y).\Eeas
We claim that this inequality can be rewritten as
\Bea P^{*,\fL}_{t_0}f(x)
& \lesssim & 
C_1(\sqrt x)\,\Ml_{M^2} \big(fe^{-\frac{y}2}\big)(x)\,+\,\frac{C_2(\sqrt x)}{x^{1/4}}\,\int_{\SRp}
|f(u)|\Phi^\fL(u).\label{P*fL}\Eea
The expression for the second term is clear from \eqref{PhifL} (after a change of variables $y^2=u$).
To handle the first term, notice that the local region now becomes $\frac {\sqrt x}2<y<M\sqrt x$,
which in particular gives $x^{-\frac 14}\sqrt y\approx1$. We also need the following lemma 
for the maximal function.

\begin{lemma}
For all  $g\geq 0$ and $x\in\SRp$,\[
\cM\Big(g(y^2)\chi_{\{\frac {\sqrt x}2<y<M\sqrt x\}}\Big)(\sqrt x) \,\lesssim\,
\cM\Big(g(u)\chi_{\{\frac x4<u<M^2x\}}\Big)(x).\]
\end{lemma}
\Proof
This follows essentially from the change of variables $y^2=u$,
\Beas 
LHS &  \leq & \sup_{r>0}\frac 1r\int_{|y-\sqrt x|<r}g(y^2)\,\chi_{\{\frac x4<y^2<M^2x\}}\,dy\\
&  = & \sup_{r>0}\frac 1r\int_{|\sqrt u-\sqrt x|<r}g(u)\,\chi_{\{\frac x4<u<M^2x\}}\,\frac{du}{2\sqrt u}.\Eeas
In this local range we have $\sqrt{u}\approx\sqrt x$, so may take the denominator outside the integral.
The local behavior also implies that\[
|\sqrt u-\sqrt x|=\Big|\int_x^u\frac{ds}{2\sqrt s}\Big|\,\approx\,\frac{|u-x|}{\sqrt x}. 
\]
Therefore, we conclude that
\[LHS\,\lesssim\, \sup_{r>0}\frac 1{r\sqrt x}\int_{|u-x|<r\sqrt x}g(u)\,\chi_{\{\frac x4<u<M^2x\}}\,du\,\lesssim\,RHS.
\]
\ProofEnd
\BR
A small variation of this proof shows that $x\in \cL_f$ implies $\sqrt x\in \cL_{Af}$, so in view of \eqref{ptwise}, the pointwise convergence in (ii) of Theorem \ref{th1SL} (for the operator $\fL$) actually holds at every Lebesgue point $x$ of $f$. 
\ER

We have thus shown \eqref{P*fL}. From here one proves Theorem \ref{th2SL} (for the operator $\fL$)
arguing once again as in $\S\ref{proofth2}$. Remark that, in view of the new constants $C_1$ and $C_2$ in \eqref{P*fL}, the weight $v=\min\{v_{1,\e},v_2\}$ must be defined with
\[
v_{1,\e}(x)=\frac{e^{-\frac{px}2}}{(1+x)^{p\nu}}\,
\mathscr{V}(x)\rho_{\e}\Big({\mathscr{V}(x)}\Big)\mand v_2(x)=\,\frac{\l{x}^{(\frac\al2+1)p-1}}{[\log(e/\l{x})]^2}\,\frac{e^{-\frac
p2x}}{(1+x)^{N}}\,,
\]
where
$\mathscr{V}(x)=\Big[\Ml_{M^2}\big(w^{\frac1{1-p}}e^{-\frac{p'x}2}\big)(x)\Big]^{1-p}$ and 
$N>1+\frac p2(|\mu+\frac12|-\frac12)$. Then, the same proof as before gives that $P^*_{t_0}:L^p(w)\to L^p(v)$ if $w\in D_p(\Phi^\fL)$. One can also establish (with a few obvious modifications) that for every given $q>p$, there is
a sufficiently small $\e$ so that $v\in D_q(\Phi^\fL)$. Once again, we may also replace this weight by 
$v^{\Phi^\fL,w}(x)$, as defined in \eqref{tvdef}.

\subsection{Results for the system $\eln$}
Remember that these functions satisfy
\Be
\eln(y)\,=\,a(y)\,\fLn(y),\quad\mbox{with } a(y)=y^{-\frac\al2}.\label{aeln}\Ee
Thus, they are eigenvectors of the differential operator
\[
f\mapsto \el f(x)=a(x)\fL[a^{-1}f](x)
\]
(with the same eigenvalues as $\fLn$) and constitute an orthonormal system in $L^2(a(y)^{-2}dy)$.
One then derives Theorems \ref{th1SL} and \ref{th2SL} for the operator $\el$, from the known
results about $\fL$, by repeating exactly the same arguments that we gave in $\S\ref{transf_psi}$.
We leave the details to the reader.

\subsection{Results for the Laguerre polynomials $\Lak$}\label{subsecLak}
This system and the corresponding operator $\SL$ in \eqref{SL} are the ones considered
in the statements of $\S\ref{intro}$, so we shall give a few more details here.
First of all, recall that $\Lak$ and $\fLn$ are linked by 
\Be
\Lak(y)\,=\,a(y)\,\fLn(y),\quad\mbox{with } a(y)=y^{-\frac\al2}e^{y/2}.\label{aLak}\Ee
Thus, the functions $\Lak$ are orthonormal in $L^2(a(y)^{-2}dy)=L^2(y^\al e^{-y}dy)$, and are also eigenvectors of the differential operator
\[
f\mapsto \SL f(x)=a(x)\fL[a^{-1}f](x),
\]
with the same eigenvalues as $\fLn$, namely $n+(\al+1+\mu)/2$. We remark that $\SL$
coincides with the operator defined in \eqref{SL} when we set $m=(\al+1+\mu)/2$. 
Thus, the heat and Poisson kernels of these two operators are related by
\[
e^{-t\SL}(x,y)\,=\,\sum_{n=0}^\infty e^{-(n+m)t}\Lak(x)\Lak(y)\,=\,
a(x)\,a(y)\,e^{-t\fL}(x,y),
\]
and
\[
P_t^{\SL}(x,y)\,=\,a(x)\,a(y)\,P_t^{\fL}(x,y).
\]
This implies the identity
\Be
P_t^\SL f(x)\,=\,\int_{\SRp}P_t^\SL(x,y)\,f(y)\,y^\al e^{-y}\,dy\,=
\,a(x)\,P^\fL_t[fy^{\frac\al2}e^{-\frac y2}](x),\label{PSLfL}\Ee
from which one deduces the validity of Theorem \ref{th1SL} for the operator $\SL$, provided 
\[
\Phi^\SL(x)\,=\,y^{\frac\al2}e^{-\frac y2}\Phi^\fL(y)\,=\,
\frac{y^\al e^{-y}}{(1+y)^{\frac{1+\al+\mu}2}[\log(e+x)]^{1+\nu}}.
\]
Note that this coincides with the function in \eqref{Phi} since we have set $m=(\al+1+\mu)/2$.
Moreover, \eqref{PSLfL} combined with \eqref{P*fL} implies the estimate
\Be P^{*,\SL}_{t_0}f(x)
\, \lesssim \,
C^\SL_1(x)\,\Ml_{M^2} \big(fe^{-{y}}\big)(x)\,+
\,C^\SL_2(x)\,\int_{\SRp}
|f(u)|\Phi^\SL(u)=\,I(x)+II(x),\label{P*SL}\Ee
with the new constants
\[
C^\SL_1(x)=(1+x)^\nu e^x\mand C^\SL_2(x)=(1+x)^{(|\mu+\frac12|-\al-\frac12)/2} e^x/\l{x}^{\al+1}.
\]We now apply the same arguments as in $\S\ref{proofth2}$ to show that, for a suitable weight $v$, we have $\|P^{*,\SL}_{t_0}f\|_{L^p(v)}\lesssim \|f\|_{L^p(w)}$, under the assumption $w\in D_p(\Phi^{\SL})$.
Indeed, to control the second term $II(x)$ we choose a weight $v_2$ such that $C_2^\SL\in L^p(v_2)$, namely
\[v_2(x)=\,\frac{\l{x}^{(\al+1)p-1}}{[\log(e/\l{x})]^2}\,\frac{e^{-px}}{(1+x)^{N}}\,,
\]
with say $N>1+p(m+|\al+\frac12|)$. It is not difficult to verify that $v_2\in D_q(\Phi^\SL)$ for all $q>p$.
To control the first term we set\[
v_{1,\e}(x)=\frac{e^{-px}}{(1+x)^{p\nu}}\,
\mathscr{V}(x)\rho_{\e}\Big({\mathscr{V}(x)}\Big)\mbox{ with }
\mathscr{V}(x)=\Big[\Ml_{M^2}\big(w^{\frac1{1-p}}e^{-p'x}\big)(x)\Big]^{1-p}.\]
That is, if $W(x)=w(x)e^{px}$, then $v_{1,\e}(x)=(1+x)^{-p\nu}e^{-px}\,V_\e(x)$
with the notation in \eqref{Ve}.
So we may quote Theorem \ref{ThVe} to obtain\[
\big\|I\big\|_{L^p(v_{1,\e})}\,\lesssim\,\Big[\int_\SRp |\Ml_{M^2} \big(fe^{-{y}}\big)(x)|^p\,V_{\e}(x)dx\Big]^{\frac1p}\,
\lesssim\,\|fe^{-{y}}\|_{L^p(W)}=\|f\|_{L^p(w)}.
\]
Again, it is not difficult to verify that for a sufficiently small $\e$ one has $v_{1,\e}\in D_q(\Phi^\SL)$,
arguing as in the last part\footnote{With the quadratic exponentials $e^{{x^2}/2}$ in \eqref{WDexp} and \eqref{aux9} replaced by linear exponetials $e^x$.} of $\S\ref{proofth2}$. Thus, Theorem \ref{th2SL} holds with $v=\min\{v_{1,\e},v_2\}$. Alternatively, with the notation in \eqref{tvdef}, one may as well choose the weight $v^{\Phi^\SL,w}(x)$.

\bline {\bf Acknowledgements:} We wish to thank Jos\'e Luis Torrea for many useful conversations around these topics at the earlier stages of this work. First and third authors partially supported by grants MTM2013-40945-P, MTM2013-42220-P and MTM2014-57838-C2-1-P from MINECO (Spain), and grants 19368/PI/14 and 19378/PI/14 from \emph{Fundaci\'on S\'eneca}, (Regi\'on de Murcia, Spain).
Second and fourth authors partially supported by grants from 
\emph{Consejo Nacional de Investigaciones Cient\'\i ficas y T\'ecnicas}
(CONICET) and Universidad Nacional del Litoral (Argentina).

\end{document}